\documentclass{amsproc}

\usepackage[arrow, curve, frame, matrix, graph, tips]{xy}
\usepackage{amssymb}
\usepackage{amsthm}
\usepackage{mtens}

\newtheorem{theorem}{Theorem}[section]
\newtheorem{lemma}[theorem]{Lemma}
\newtheorem{proposition}[theorem]{Proposition}
\newtheorem{corollary}[theorem]{Corollary}

\theoremstyle{definition}
\newtheorem{definition}[theorem]{Definition}
\newtheorem{example}[theorem]{Example}
\newtheorem{examples}[theorem]{Examples}
\newtheorem{non-example}[theorem]{Non-Example}

\theoremstyle{remark}
\newtheorem{remark}[theorem]{Remark}

\begin{document}

\title{Free Products of Higher Operad Algebras}

\author{Mark Weber}
\address{Max Planck Institute for Mathematics, Bonn}
\email{weber@pps.jussieu.fr}
\thanks{}
\maketitle
\maketitle
\begin{abstract}
One of the open problems in higher category theory is the systematic construction of the higher dimensional analogues of the Gray tensor product of $2$-categories. In this paper we continue the developments of \cite{EnHopI} and \cite{EnHopII} by understanding the natural generalisations of Gray's little brother, the funny tensor product of categories. In fact we exhibit for any higher categorical structure definable by an $n$-operad in the sense of Batanin \cite{Bat98}, an analogous tensor product which forms a symmetric monoidal closed structure on the category of algebras of the operad.
\end{abstract}
\tableofcontents

\section{Introduction}\label{sec:introduction}
Strict $n$-categories have an easy inductive definition, with a strict $(n{+}1)$-category being a category enriched in the category of strict $n$-categories via its \emph{cartesian} product, but unfortunately these structures are too strict for intended applications in homotopy theory and geometry. In dimension 2 as is well-known, there is no real problem because any bicategory is biequivalent to a strict 2-category. However in dimension 3 the strictest structure one can replace an arbitrary weak $3$-category with -- and not lose vital information -- is a Gray category, which is a category enriched in $\Enrich 2$ using the Gray tensor product of 2-categories instead of its cartesian product \cite{GPS95}. This leads naturally to the idea of trying to define what the higher dimensional analogues of the Gray tensor product are, in order to set up a similar inductive definition as for strict $n$-categories, but to capture the strictest structure one can replace an arbitrary weak $n$-category with and \emph{not} lose vital information.

Ignoring the 2-cells in the definition of the Gray tensor product one has a canonical tensor product of categories, which has been called the \emph{funny tensor product}. It is different from the cartesian product of categories. In this article we explain that for \emph{any} higher categorical structure definable by an $n$-operad $A$ in the sense of \cite{Bat98}, one has an analogous tensor product giving the category of $A$-algebras a symmetric monoidal closed structure. We call these generalisations of the funny tensor product ``free products'' for reasons that will become clear.

In \cite{EnHopII} the foundations of higher category theory in the globular style were given a major overhaul. The resulting theory is expressed as the interplay between lax monoidal structures on a category $V$ and monads on the category $\ca GV$ of graphs enriched in $V$. This perspective enables one to understand many issues in higher category theory without looking at the combinatorics that might, at first glance, appear to make things difficult. Our description and analysis of the generalisations of the funny tensor product is an illustration of this. Thus the present paper is written within the framework of established by \cite{EnHopII}, and so we shall use the terminology and notations of that paper here without further comment.

We recall the funny tensor product of categories in section(\ref{sec:funny}), focussing in particular on those aspects which we have been able to generalise to all higher categorical structures in place of $\Cat$. Then in section(\ref{sec:monoidal-monads}) we review the theory of symmetric monoidal monads from the point of view of multicategories. None of this section is new. The theory of symmetric monoidal monads is originally due to Anders Kock \cite{KockMM1} \cite{KockMM2} \cite{KockMM3}, and in the setting of Lawvere theories to Fred Linton \cite{LintonCommTh}. However, our use of multicategorical notions, especially representablility in the sense of Hermida \cite{Her00} and closedness in the sense of Manzyuk \cite{Manzyuk}, to help understand monoidal monads does appear to be original. It is this synthesis of the theories of symmetric monoidal monads and of multicategories that provides the most convenient framework within which to construct tensor products from monads and operads.

A notion of multimap of enriched graph is provided by a general construction in section(\ref{ssec:obj-F}) which applies to any category equipped with a functor into $\Set$. The 2-functoriality of this construction described in section(\ref{ssec:2functoriality-F}) and its compatibility with monad theory described in section(\ref{ssec:monads-from-F}), enables one to give any monad defined on $\ca GV$ over $\Set$ a canonical symmetric monoidal structure. Applying then our formulation of the theory of symmetric monoidal monads, we are able to exhibit our higher dimensional analogues of the funny tensor product in section(\ref{ssec:fp-monad-algebras}). Some of the formal properties that they enjoy are exhibited in sections(\ref{ssec:free-vs-cartesian}) and (\ref{ssec:pushout-formula}).

In particular in section(\ref{ssec:free-vs-cartesian}) we isolate a condition on a monad so that one obtains a canonical identity on objects comparison between the free product of $T$-algebras and its cartesian product. This important because, as explained to the author by John Bourke, the Gray tensor product for 2-categories may be obtained by factorising this map, in the case where $T$ is the monad on $\ca G^2\Set$ for 2-categories, using the bijective on objects fully faithful factorisation for $\Cat$ on the homs. In this way the coherence data for the Gray tensor product is determined by what we know more generally from our theoretical framework and the magic of orthogonal factorisation systems. For this reason it seems that a complete understanding of the higher dimensional analogues of the Gray tensor product is within reach.

Given an $n$-operad $A$ one may consider categories enriched in the algebras of $A$ for the free product. We call such structures \emph{sesqui-$A$-algebras}. In section(\ref{ssec:monads-sesqui-algebras}) we give an explicit description of the monad whose algebras are sesqui-$A$-algebras. This depends on a general result, called the ``multitensor dropping theorem'', given in section(\ref{ssec:dropping-theorem}). Then in section(\ref{ssec:operads-sesqui-algebras}) we explain why this monad is part of an $(n{+}1)$-operad. These results will form part of the inductive machine for semi-strict $n$-categories that we hope to uncover in future work.

\section{The funny tensor product of categories}\label{sec:funny}
The category $\Cat$ is a cartesian closed category. By
\begin{theorem}
\emph{Foltz, Kelly and Lair}\cite{FKL}
Up to isomorphism there are exactly two biclosed monoidal structures on $\Cat$, and both are symmetric.
\end{theorem}
\noindent $\Cat$ has another symmetric monoidal closed structure. The corresponding tensor product has been called the \emph{funny tensor product}. It is related to the cartesian product by identity-on-objects functors
\[ \kappa_{A,B} : A \tensor B \rightarrow A \times B \]
which are natural in $A$ and $B$. Thus the objects of $A \tensor B$ are pairs $(a,b)$ with $a \in A$ and $b \in B$. A generators and relations description of the morphisms of $A \tensor B$ is as follows. They are generated by
\[ \begin{array}{lccr} {(a,\beta) : (a,b_1) \rightarrow (a,b_2)} &&& {(\alpha,b) : (a_1,b) \rightarrow (a_2,b)} \end{array} \]
where $a$ and $\alpha:a_1{\rightarrow}a_2$ are in $A$, and $b$ and $\beta:b_1{\rightarrow}b_2$ are in $B$. The relations are obtained by remembering composition in $A$ and $B$, that is
\[ \begin{array}{lcccr} {(\alpha_1\alpha_2,b) = (\alpha_1,b)(\alpha_2,b)} && {(1_a,b)=1_{(a,b)}=(a,1_b)} && {(a,\beta_1\beta_2) = (a,\beta_1)(a,\beta_2).} \end{array} \]
So given maps $\alpha$ and $\beta$ one has the square
\[ \xygraph{!{0;(3,0):(0,.333)::} {(a_1,b_1)}="tl" [r] {(a_1,b_2)}="tr" [d] {(a_2,b_2)}="br" [l] {(a_2,b_1)}="bl" "tl"(:"tr"^-{(a_1,\beta)}:"br"^{(\alpha,b_2)},:"bl"_{(\alpha,b_1)}:"br"_-{(a_2,\beta)})} \]
which by contrast with the cartesian product does not commute in general. In other words this square in $A \tensor B$ has not one but \emph{two diagonals}, the composites $(a_2,\beta)(\alpha,b_1)$ and $(\alpha,b_2)(a_1,\beta)$ which are identified by the functor $\kappa_{A,B}$. An explicit description of the Gray tensor product of $2$-categories proceeds in the same way as for $\tensor$ on objects and arrows, with the key feature in dimension 2 being an isomorphism between these two diagonals.

While $\tensor$ acts on objects like cartesian product, on morphisms it behaves more like a \emph{co}product. In particular taking $A$ and $B$ to be monoids, that is to say one object categories, then $A \tensor B$ is of course also a monoid, and is in fact the coproduct, in the category of monoids, of $A$ and $B$. It is standard terminology from algebra to refer to the coproduct in the category of monoids as the \emph{free product} of monoids, and so after this section we shall adopt this terminology, referring to $\tensor$ and its generalisations as ``free products''.

For general categories this coproduct-like behaviour for the arrows of $A \tensor B$ is expressed by the pushout
\[ \xygraph{!{0;(2,0):(0,.5)::} {A_0{\times}B_0}="tl" [r] {A_0{\times}B}="tr" [d] {A{\tensor}B}="br" [l] {A{\times}B_0}="bl" "tl"(:"tr"^-{\id{\times}i_B}:@{.>}"br",:"bl"_{i_A{\times}\id}:@{.>}"br") "br" [u(.3)l(.3)] (:@{-}[r(.15)],:@{-}[d(.15)])} \]
in $\Cat$ in which $A_0$ (resp. $B_0$) denote the underlying discrete category of $A$ (resp. $B$), and $i_A$ (resp. $i_B$) are the inclusions. We shall call this the \emph{pushout formula} for the funny tensor product $A \tensor B$.

The internal hom $[A,B]$ of $A$ and $B$ corresponding to the funny tensor product, has functors $A \rightarrow B$ as objects just as in the cartesian case, but mere \emph{transformations} as morphisms. A transformation $\phi:f \rightarrow g$ between functors $A \rightarrow B$ consists of a morphism $\phi_a:fa \rightarrow ga$ for each $a \in A$. These components $\phi_a$ are not required to satisfy the usual naturality condition as in the cartesian case.

\section{The theory of monoidal monads and multicategories}\label{sec:monoidal-monads}

\subsection{Motivation}\label{ssec:motivation-sm-review}
In this section we will describe what for us is the ``abstract categorical theory of tensor products which arise from universal properties''. Throughout this section the reader is invited to keep in mind the basic guiding example of the symmetric monoidal closed structure on the category $\RMod$ of modules over a commutative ring $R$.

It was Claudio Hermida who expressed in \cite{Her00} how one may regard this monoidal structure as arising from the \emph{multicategory} of $R$-modules and $R$-multilinear maps between them. So for Hermida this multicategory is the more fundamental object, with the monoidal structure on $\RMod$ just an expression of its \emph{representability} in his sense. This notion captures abstractly the idea that the tensor product of $R$-modules, by definition, classifies $R$-multilinearity.

A different perspective comes from the theory of symmetric monoidal monads{\footnotemark{\footnotetext{Actually Kock orginally called them \emph{commutative monads}.}}} of Anders Kock \cite{KockMM1} \cite{KockMM2} \cite{KockMM3}. By definition a symmetric monoidal monad is a monad in the 2-category $\SYMMON$ of symmetric monoidal categories, lax symmetric monoidal functors and monoidal natural transformations. Thus the data of a symmetric monoidal monad is
\[ (V,\tensor,T,\phi,\eta,\mu) \]
where $(V,\tensor)$ is a symmetric monoidal category, $(T,\phi)$ is a lax symmetric monoidal functor, the components of $\phi$ look like
\[ \begin{array}{c} {\phi_{X_i} : \Tensor\limits_i TX_i \rightarrow T\Tensor\limits_iX_i,} \end{array} \]
and $\eta$ and $\mu$ are the unit and multiplication of the monad and are monoidal transformations. The precise example which relates to the present discussion is
\begin{example}\label{ex:KockMnd}
Regard a commutative ring $R$ as a monad on $\Set$, where for $X \in \Set$, $R(X)$ is the set of formal $R$-linear combinations of elements of $X$. The unit is given by $x \mapsto 1 \cdot x$, and the multiplication of the monad is given by the evident substitution of formal linear combinations. An algebra of $R$ is simply an $R$-module. Given a finite sequence of sets $X_1,...,X_n$ one has a function
${\prod\limits_i RX_i \rightarrow R\prod\limits_i X_i}$
defined by
\[ \begin{array}{lcr} {(\sum\limits_{j_i=1}^{m_i} \lambda_{ij_i}x_{ij_i} \,\, : \,\, 1{\leq}i{\leq}n)} & {\mapsto} & 
{\sum\limits_{j_1,...,j_n} (\prod\limits_i \lambda_{ij_i}) (x_{1j_1},...,x_{nj_n})} \end{array} \]
which provides the monoidal functor coherences.
\end{example}
\noindent For Kock the above monad $R$ is the fundamental combinatorial object from which the symmetric monoidal closed structure of $\RMod$ may be obtained.

In this section we review the theory of symmetric monoidal monads from the point of view of multicategories. Given a symmetric monoidal monad $(V,\tensor,T,\phi,\eta,\mu)$, we will first see that $V^T$ underlies a symmetric multicategory. Then under some hypotheses this multicategory will be seen to be closed in the sense of \cite{Manzyuk}. Closed multicategories are more easily exhibited as representable, as observed in section(\ref{ssec:multicat-review}), and we use this to recover the tensor product on $V^T$ by exhibiting the corresponding multicategory of algebras as representable.

So for us the main results of symmetric monoidal monad theory are expressed as the closedness and representablilty of the symmetric multicategory of algebras of the original symmetric monoidal monad. Thus in our treatment, both Kock's and Hermida's perspectives on ``how the symmetric monoidal closed structure of $\RMod$ arises'' are placed on an equal footing.

\subsection{Notation and terminology} Depending on what is most convenient in a given situation
\[ \begin{array}{lcccr} {(X_i)_{1{\leq}i{\leq}n}} && {(X_i)_i} && X \end{array} \]
are different notations we shall use for the same thing, namely a sequence of sets $(X_1,...,X_n)$. Similarly a typical element $(x_1,...,x_n)$ of the cartesian product of these sets has the alternative notations
\[ \begin{array}{lcccr} {(x_i)_{1{\leq}i{\leq}n}} && {(x_i)_i} && {x.} \end{array} \]
For example, we may speak of multimaps $f:x{\rightarrow}y$ in a given multicategory $Z$, where $x$ is a \emph{sequence} of objects and $y$ is a single object from $Z$, and the set of such may be denoted as $Z(x,y)$. When the sequence $x$ has length $1$, we say that the multimap $f$ is \emph{linear}, and we denote by $\lin(Z)$ the category of objects and linear maps in $Z$, this being the object part of a 2-functor
\[ \lin : \SYMMULT \rightarrow \CAT \]
out of the 2-category of $\SYMMULT$ of symmetric multicategories.
We denote elements of sequences of sequences of sets in a similar manner consistent with our tensor product notation. For example given multimaps
\[ \begin{array}{lccr} {f_i:(x_{i1},...,x_{in_1}) \rightarrow y_i} &&& {g:(y_1,...,y_k) \rightarrow z} \end{array} \]
in a multicategory, one usually denotes their composite as
\[ g(f_1,...,f_k) : (x_{11},...,x_{1n_1},......,x_{k1},...,x_{kn_k}) \rightarrow z \]
but we shall sometimes use the notation
\[ g(f_i)_i : (x_{ij})_{ij} \rightarrow z \]
for the same thing. By way of illustration let us recall the definition of the objects, arrows and 2-cells of the 2-category $\SYMMULT$. A symmetric multicategory $X$ consists of the following data:
\begin{enumerate}
\item a set $X_0$ whose elements are called \emph{objects}.
\item for each sequence $x$ and element $y$ from $X_0$, a set $X(x,y)$ whose elements are called \emph{multimaps from $x$ to $y$}.
\item for all $f:x{\rightarrow}y$ and $\sigma \in \Sigma_{l(x)}$, a multimap $f\sigma:x\sigma{\rightarrow}y$.
\item for all $x$ in $X_0$, a distinguished multimap $1_x:(x){\rightarrow}x$ called the \emph{identity} for $x$.
\item given $f_i:(x_{ij})_j{\rightarrow}y_i$ and $g:(y_i)_i{\rightarrow}z$, another multimap $g(f_i)_i: (x_{ij})_{ij} \rightarrow z$ called the \emph{composite} of $g$ and the $f_i$.
\end{enumerate}
This data must satisfy the following axioms:
\begin{enumerate}
\item \emph{unit law of symmetric group actions}: for all $f:x{\rightarrow}y$, $f1_{l(x)}=f$.
\item \emph{associativity of symmetric group actions}: for all $f:x{\rightarrow}y$ and $\sigma,\tau$ in $\Sigma_{l(x)}$, $(f\sigma)\tau = f(\sigma\tau)$.
\item \emph{unit laws of composition}: for all $f:x{\rightarrow}y$, $1_y(f)=f=f(1_{x_i})_i$.
\item \emph{associativity of composition}: given
\[ \begin{array}{lcccr} {f_{ij}:(x_{ijk})_k{\rightarrow}y_{ij}} && {g_i:(y_{ij})_j{\rightarrow}z_i} && {h:(z_i)_i{\rightarrow}w} \end{array} \]
one has $(h(g_i)_i)(f_{ij})_{ij} = h(g_i(f_{ij})_j)_{i}$.
\item \emph{equivariance}: given $f_i:(x_{ij})_j{\rightarrow}y_i$ and $g:(y_i)_i{\rightarrow}z$, $\sigma \in \Sigma_{l(y_i)_i}$ and $\tau_i \in \Sigma_{l(x_{ij})_j}$, one has $(g(f_i)_i)\sigma(\tau_i)_i = (g\sigma)(f_{\sigma_i}\tau_{\sigma{i}})_i$.
\end{enumerate}
In the absence of the symmetric group actions and the equivariance of composition, one has the definition of a \emph{multicategory}. Let $X$ and $Y$ be symmetric multicategories. A \emph{symmetric multifunctor} $F:X{\rightarrow}Y$ consists of the following data:
\begin{enumerate}
\item  A function $F_0:X_0{\rightarrow}Y_0$. For a sequence $x$ from $X$, we abuse notation and write $Fx$ for the sequence $(F_0x_1,...,F_0x_n)$, and for an object $y$ of $X$, we write $Fy$ for the object $F_0y$ of $Y$.
\item  For each sequence $x$ from $X$ and element $y \in X_0$, a function \[ F_{x,y}:X(x,y){\rightarrow}Y(Fx,Fy). \] For a given multimap $f:x{\rightarrow}y$ in $X$, we denote by $Ff$ the multimap $F_{x,y}(f)$ in $Y$.
\end{enumerate}
This data must satisfy the following axioms:
\begin{enumerate}
\item  For all $f:x{\rightarrow}y$ in $X$ and $\sigma \in \Sigma_{lx}$, $F(f\sigma)=F(f)\sigma$.
\item  For all $x \in X_0$, $F(1_x)=1_{Fx}$.
\item  For all multimaps $f_i:(x_{ij})_j{\rightarrow}y_i$ and $g:(y_i)_i{\rightarrow}z$ in $X$, $F(g(f_i)_i) = Fg(Ff_i)_i$.
\end{enumerate}
If $X$ and $Y$ are mere multicategories and the first condition is ignored, then one has the definition of \emph{multifunctor}. Let $F$ and $G$ be symmetric multifunctors $X{\rightarrow}Y$. Then a multinatural transformation $\phi:F{\rightarrow}G$ consists of a multimap $\phi_x:F(x){\rightarrow}Gx$ for each $x \in X_0$. This data is required to satisfy the condition, called \emph{multinaturality}, that for all multimaps $f:x{\rightarrow}y$ in $X$, one has $Gf(\phi_{x_i})_i=\phi_y(Ff)$.

Returning to our general element $(x_i)_i$ of $\prod\limits_iX_i$, if $z$ is an element or sequence of elements of $X_k$ for some $k$, then we denote by $x|_kz$ the new sequence of elements obtained by replacing $x_k$ with $z$.  

\subsection{Closed and representable symmetric multicategories}\label{ssec:multicat-review}
Let us recall the forgetful 2-functor
\[ U : \SYMMON \rightarrow \SYMMULT{.} \]
Given a symmetric monoidal category $(V,\tensor)$, the symmetric multicategory $UV$ has the same objects as $V$ and homs given by
\[ \begin{array}{c} {UV(x,y) = V(\Tensor\limits_ix_i,y)} \end{array} \]
and the rest of the structure is given in the obvious way using the symmetric monoidal category structure. Given a symmetric lax monoidal functor $(F,\phi):V{\rightarrow}W$, one defines $UF$ to have the same object map as $F$ and hom maps given by
\[ \begin{array}{lcr} {f:\Tensor\limits_ix_i{\rightarrow}y} & \mapsto &
{\xygraph{!{0;(1.5,0):} {\Tensor\limits_iFx_i}="l" [r] {F(\Tensor\limits_ix_i)}="m" [r] {Fy}="r" "l":@<.5ex>"m"^-{\phi_{x_i}}:@<.5ex>"r"^-{Ff}}.} \end{array} \]
Given a monoidal natural transformation $\psi:(F,\phi){\rightarrow}(G,\gamma)$, the induced multinatural transformation $U\psi$ has components which can be identified with those of $\psi$ once the linear maps of $UW$ are identified with the morphisms of $W$. The verification that these assignments form a 2-functor is routine and left to the reader. For us the key facts about $U$ are summarised by
\begin{proposition}\label{prop:smultcat->smoncat}
Let $X$ be a symmetric multicategory and $U$ the 2-functor just described.
\begin{enumerate}
\item $X \iso UV$ for some symmetric monoidal category $V$ iff $X$ is representable in the sense of Hermida \cite{Her00}.\label{U-rep}
\item $U$ is 2-fully-faithful.\label{U-2ff}
\end{enumerate}
\end{proposition}
\noindent This result is a mild reformulation of some aspects of the theory of multicategories already in the literature, especially in \cite{Her00}. Let us briefly recall some of this theory, and in so doing, explain why this result is true.

Let $X$ be a multicategory. A multimap $f:x{\rightarrow}y$ in $X$ is \emph{universal} when for all $z \in X$ the function
\[ (-) \comp f : X(y,z) \rightarrow X(x,z) \]
given by composition with $f$ is a bijection. In other words any multimap out of $x$ may be identified with a unique linear map out of $y$. The example to keep in mind is the multicategory of vector spaces and multilinear maps. In this case the universality of $f$ expresses that the vector space $y$ is exactly the tensor product of the $x_i$. Thus one might be tempted to guess that a monoidal category is really just a multicategory in which every tuple of objects admits such a universal multimap. While this is a good first guess, it is not quite correct in general, and it was Hermida who understood in \cite{Her00} how to express this idea correctly.

Hermida's key insight is that one must strengthen the notion of universal map. In order to express this idea it is necessary to recall placed composition in a multicategory. Given a multimap $g:y{\rightarrow}z$ and $f:x{\rightarrow}y_k$ for some $k$, one defines their \emph{placed composite} $g \comp_k f : y|_kx \rightarrow z$ as the composite
\[ g(1_{y_1},...,1_{y_{k-1}},f,1_{y_{k{+}1}},...,1_{y_{n}}) \]
in $X$. Clearly one can recapture general composition from placed composition. A multimap $f:x{\rightarrow}y$ in $X$ is \emph{strongly universal} when for all sequences $y'$ and $k$ such that $y'_k=y$, and objects $z$ of $X$, the function
\[ (-) \comp_k f : X(y',z) \rightarrow X(y'|_kx,z) \]
given by placed composition with $f$ is a bijection. The multicategory $X$ is \emph{representable} when for every tuple of objects $(x_i)_i$ of $X$, there exists a strongly universal multimap out of it. In order to bring out the relation with monoidal categories one may denote a choice of such a strongly universal multimap as
\[ \begin{array}{c} {(x_i)_i \rightarrow \Tensor\limits_ix_i} \end{array} \]
to emphasize that from such a choice one obtains a monoidal structure on the category of objects and linear maps of $X$, and the non-symmetric analogue of $U$ sends this monoidal category to $X$. It is not very hard to adapt these insights to the symmetric case. A symmetric multicategory is representable when it is representable as a mere multicategory, and when $X$ is a representable symmetric multicategory with a given choice of universal maps, the symmetries on $X$ give rise in the evident way to symmetries for the induced tensor product on the linear part of $X$.

As for the 2-fully-faithfulness, Hermida showed that monoidal categories are the adjoint pseudo-algebras of a lax idempotent 2-monad on the 2-category of multicategories. This too is easily adapted to the symmetric case, and so one may regard $U$ as the forgetful 2-functor from the 2-category of adjoint pseudo-algebras, lax maps and algebra 2-cells, for a lax idempotent 2-monad on $\SYMMULT$. By the general theory of such 2-monads, see for instance \cite{KL97}, such forgetful 2-functors are always 2-fully-faithful.

The distinction between universal and strongly universal multimaps disappears when one works with closed symmetric multicategories, as we shall see in the next lemma. Closed multicategories, and their relation to closed categories, were discussed in \cite{Manzyuk}. We recall this notion in the symmetric context. Let $X$ be a symmetric multicategory. Then $X$ is \emph{closed} when for all objects $x$ and $y$ in X, there is an object $[x,y]$ in $X$ together with a ``right evaluation'' multimap
\[ \rev_{x,y} : ([x,y],x) \rightarrow y \]
such that for all sequences $(z_1,...,z_n)$ of objects of $X$, the function
\[ \rev_{x,y}(-,1_x) : X((z_1,...,z_n),[x,y]) \rightarrow X((z_1,...,z_n,x),y) \]
is a bijection. This notion is not as one-sided as it may first appear because of the presence of symmetries in $X$, enabling one to define ``left evaluation'' multimaps
\[ \lev_{x,y} : (x,[x,y]) \rightarrow y \]
with the property that for all sequences $(z_1,...,z_n)$ of objects of $X$, the function
\[ \lev_{x,y}(1_x,-) : X((z_1,...,z_n),[x,y]) \rightarrow X((x,z_1,...,z_n),y) \]
is a bijection.
\begin{remark}
By definition a representable symmetric multicategory is closed iff the corresponding symmetric monoidal category is closed in the usual sense.
\end{remark}
\begin{lemma}\label{lem:closed->representable}
Let $X$ be a closed symmetric multicategory and $f:(x_1,...,x_p){\rightarrow}y$ be a multimap therein. If $f$ is universal then $f$ is in fact strongly universal.
\end{lemma}
\begin{proof}
For $m,n \in \N$ and sequences of objects $(z_1,...,z_{m+n})$ of $X$, we must show that the function
\[ X(z_1,...,z_m,y,z_{m{+}1},...,z_{m{+}n},w) \rightarrow X(z_1,...,z_m,x_1,...,x_p,z_{m{+}1},...,z_{m{+}n},w) \]
which we denote by $\phi_{m,n,w}$, given by
\[ \phi_{m,n,w}(g) = g(1_{z_1},...,1_{z_m},f,1_{z_{m{+}1}},...,1_{z_{m{+}n}}) \]
is bijective for all objects $w$ of $X$. Since $X$ is symmetric it suffices to consider the case $n=0$. We proceed by induction on $m$. The base case $m=0$ follows since $f$ is universal by assumption. For the inductive step one exhibits the following string of bijections:
\[ \xygraph{!{0;(1,0):(0,.5)::} {(z_1,...,z_{m{+}1},y) \rightarrow w}="t1" [d]
{(z_2,...,z_{m{+}1},y) \rightarrow [z_1,w]}="t2" [d]
{(z_2,...,z_{m{+}1},x_1,...,x_p) \rightarrow [z_1,w]}="t3" [d]
{(z_1,...,z_{m{+}1},x_1,...,x_p) \rightarrow w}="t4"
"t1":@{}"t2"|{}="d1":@{}"t3"|{}="d2":@{}"t4"|{}="d3"
"d1" [l(3)] :@{-}[r(6)] [r] {\textnormal{(closedness)}}
"d2" [l(3)] :@{-}[r(6)] [r] {(\phi_{m,0,z_1\backslash{w}})}
"d3" [l(3)] :@{-}[r(6)] [r] {\textnormal{(closedness)}}} \]
the composite assignment being that of $\phi_{m{+}1,0,w}$ by the multinaturality of composition in $X$.
\end{proof}
\noindent Thus it is easier to show that a symmetric multicategory is representable when one knows in advance that the symmetric multicategory in question is closed. For the non-symmetric version of this result, one must work with the 2-sided analogue of the notion discussed in \cite{Manzyuk}, that is, what one would call \emph{biclosed multicategories}. We shall not pursue this any further here.

In \cite{Manzyuk} Manzyuk realised that the only obstruction to obtaining a closed category from a closed multicategory is the existence of a unit, and so he made
\begin{definition}\cite{Manzyuk}
A \emph{closed symmetric multicategory with unit} is a closed symmetric multicategory $X$ together with a universal multimap $u : () \rightarrow e$ therein. 
\end{definition}
\noindent in the non-symmetric context. By lemma(\ref{lem:closed->representable}) such a map $u$ is strongly universal.

Recall that a closed structure on a category $V$ consists of an object $I$ called the unit, a functor
\[ [-,-] : \op V \times V \rightarrow V \]
called the hom, and morphisms
\[ \begin{array}{lcccr} {i_A : A \rightarrow [I,A]} && {j_A : I \rightarrow [A,A]} && {L^A_{B,C} : [B,C] \rightarrow [[A,B],[A,C]]} \end{array} \]
natural in their arguments, such that the $i_A$ are isomorphisms and
\[ \begin{array}{lccr} {\xybox{\xygraph{{I}="l" [r(2)] {[B,B]}="r" [dl] {[[A,B],[A,B]]}="b" "l"(:"r"^-{j}:"b"^{L},:"b"_{j})}}} &&&
{\xybox{\xygraph{{[A,B]}="l" [r(2)] {[I,[A,B]]}="r" [dl] {[[A,A],[A,B]]}="b" "l"(:"r"^-{i}:@{<-}"b"^{[j,\id]},:"b"_{L})}}} \end{array} \]
and
\[ \xygraph{{\xybox{\xygraph{!{0;(1.25,0):(0,.666)::} {[C,D]}="tl" [r(2)] {[[B,C],[B,D]]}="tr" [d] {[[B,C],[[A,B],[A,D]]]}="mr" [dl] {[[[A,B],[A,C]],[[A,B],[A,D]]]}="b" [ul] {[[A,C],[A,D]]}="ml"
"tl"(:"tr"^-{L}:"mr"^{[\id,L]}:@{<-}"b"^(.4){[L,\id]},:"ml"_{L}:"b"_(.4){L})}}}
[r(4.5)]
{\xybox{\xygraph{!{0;(.8,0):(0,1.25)::} {[A,B]}="l" [r(2)] {[A,[I,B]]}="r" [dl] {[[I,A],[I,B]]}="b" "l"(:"r"^-{[\id,i]},:"b"_(.4){L}:"r"_(.6){[i,\id]})}}}
} \]
are commutative.

Given a closed symmetric multicategory $X$ with unit $u:(){\rightarrow}e$ one exhibits a closed structure on $\lin(X)$ as follows. The unit is $e$ and the hom is just the hom of $X$. For $a \in X$ the linear morphisms $i_a$ and $j_a$ are defined uniquely by the requirement that
\[ \begin{array}{lccr}
{\xybox{\xygraph{!{0;(1.5,0):(0,.666)::} {a}="tl" [r] {(a,e)}="tr" [d] {([e,a],e)}="br" [l] {a}="bl" "tl"(:"bl"_{\id},:"tr"^-{(\id,u)}:"br"^{(i_a,\id)}:"bl"^-{\rev})}}} &&
{\xybox{\xygraph{!{0;(1.5,0):(0,.666)::} {a}="tl" [r] {(e,a)}="tr" [d] {([a,a],a)}="br" [l] {a}="bl" "tl"(:"bl"_{\id},:"tr"^-{(u,\id)}:"br"^{(j_a,\id)}:"bl"^-{\rev})}}}
\end{array} \]
be commutative, using the universal property of $\rev$ and the strong universality of $u$. For $a,b,c \in C$ the linear morphism $L^a_{b,c}$ is defined uniquely by the requirement that
\[ \xygraph{!{0;(1.5,0):(0,.666)::} {([b,c],[a,b],a)}="tl" [r(2)] {([b,c],b)}="tm" [r(2)] {c}="tr" [dl] {([a,c],a)}="br" [l(2)] {([[a,b],[a,c]],[a,b],a)}="bl"
"tl"(:"tm"^-{(\rev,\id)}:"tr"^-{\rev},:"bl"_(.4){(L^a_{b,c},\id,\id)}:"br"_-{(\rev,\id)}:"tr"_{\rev})} \]
is commutative, using the universal property of $\rev$ (twice). The result
\begin{proposition}\cite{Manzyuk}\label{prop:closedmcat->closedcat}
Let $X$ be a closed symmetric multicategory with unit $u:(){\rightarrow}e$. Then the data \[ (e,[-,-],i,j,L) \] defines a closed category structure on $\lin(V)$.
\end{proposition}
\noindent is a special case of \cite{Manzyuk} proposition(4.3).

\subsection{Symmetric monoidal monads}\label{ssec:smonmnd} As mentioned in section(\ref{ssec:motivation-sm-review}) the data of a symmetric monoidal monad is \[ (V,\tensor,T,\phi,\eta,\mu) \] where $(V,\tensor)$ is a symmetric monoidal category, $(T,\phi)$ is a lax symmetric monoidal functor, and $\eta$ and $\mu$ are the unit and multiplication of the monad and are monoidal natural transformations. By proposition(\ref{prop:smultcat->smoncat}), applying $U$ gives a bijection between symmetric monoidal monads on $(V,\tensor)$ and monads in $\SYMMULT$ on $UV$.

As a 2-category $\SYMMULT$ has all limits and colimits \cite{WeissThesis} and thus in particular it has Eilenberg-Moore objects. However one may readily exhibit these directly. Let $(X,T,\eta,\mu)$ be a monad in $\SYMMULT$: $X$ is the symmetric multicategory on which it acts, $T$ is the underlying multiendofunctor, and $\eta$ and $\mu$ are the unit and multiplication multinatural transformations. An algebra of $T$ is defined to be an algebra of the monad $\lin(T)$: that is a pair $(x,a)$ where $x$ is the underlying object in $X$ and $a:Tx{\rightarrow}x$ the action satisfying the usual axioms. Let \[ (x_1,a_1),...,(x_n,a_n),(y,b) \] be $T$-algebras and $f:(x_i)_i{\rightarrow}y$ be in $X$. Then $f$ is a \emph{multi-$T$-algebra morphism} when it satisfies the equation $bT(f)=f(a_i)_i$. Equivalently one can express this last condition on $f$ one variable at a time. That is given $1 \leq j \leq n$, $f$ is said to be a \emph{$T$-algebra morphism in the $j$-th variable} when it satisfies
\[ f(1_{x_i}|_ja_j)_i = bT(f)(\eta_{x_i}|_j1_{Tx_j})_i, \]
the domain of this composite multimap being the sequence $(x_i|_jTx_j)_i$. By substituting $(\eta_{x_i}|_j1_{Tx_j})_i$ into $bT(f)=f(a_i)_i$, one sees that $f$ is a multi-$T$-algebra morphism implies that it is a $T$-algebra morphism in the $j$-th variable. Conversely, if $f$ is a $T$-algebra morphism in each variable, then a straight forward inductive argument may be given to show that $f$ is in fact a multi-$T$-algebra morphism{\footnotemark{\footnotetext{We emphasize that the base case for the induction is $n=0$. That is the statement: $f$ is a $T$-algebra morphism in each variable iff it is a multi-$T$-algebra morphism, is true when $n=0$. This makes it a vacuous condition in this case which may seem strange. The reason for this is that the \emph{multi}-naturality condition of $\eta$, in fact of any multinatural transformation, is in particular valid with respect to nullary multimaps, and so the multi-$T$-algebra morphism condition in this case amounts to the unit law for $(y,b)$.}}}.

One defines the multicategory $X^T$ of algebras of $T$ with $T$-algebras as objects and multi-$T$-algebra morphisms as multimaps. With the evident composition and symmetric group actions inherited from $X$, one may easily verify that $X^T$ is indeed a symmetric multicategory and that one has a forgetful symmetric multifunctor
\[ \begin{array}{lccr} {U^T : X^T \rightarrow X} &&& {(x,a) \mapsto x} \end{array} \]
with object map as indicated in the previous display. Moreover one has a multinatural transformation
\[ \xygraph{!{0;(2,0):(0,.3)::} {X^T}="l" [ur] {X}="t" [d(2)] {X}="b" "l"(:"t"^{U^T}:"b"^{T},:"b"_{U^T}|{}="precodtau" "t":@{}"precodtau"|(.5){}="dtau"|(.75){}="ctau")
"dtau":@{=>}"ctau"^{\tau}} \]
whose component at $(x,a)$ is just $a$. The verification of the following result is a mild variation of the analogous easy result for monads in $\CAT$, and is left to the reader.
\begin{proposition}\label{prop:EM-SYMMULT}
Let $(X,T,\eta,\mu)$ be a monad in the 2-category $\SYMMULT$. Then its Eilenberg-Moore object in $\SYMMULT$ is given by $(U^T,\tau)$.
\end{proposition}
\noindent Thus by applying $U$ to a symmetric monoidal monad $(V,\tensor,T,\phi,\eta,\mu)$ and taking Eilenberg-Moore objects in $\SYMMULT$, one finds that the category of algebras $V^T$ is in fact the category of linear maps of a symmetric multicategory. In fact \[ V^T = \lin((UV)^{(UT)}). \]
Unpacking the notion of multimap of $T$-algebras so obtained, one finds that a $T$-algebra multimap $(X_i,x_i)_i \rightarrow (Y,y)$ consists of a morphism
\[ \begin{array}{c} {f : \Tensor\limits_i X_i \rightarrow Y} \end{array} \]
in $V$ such that
\[ \xygraph{{\Tensor\limits_i TX_i}="tl" [r(2)] {T\Tensor\limits_iX_i}="tr" [r(2)] {TY}="mr" [dl] {Y}="b" [l(2)] {\Tensor\limits_iX_i}="ml"
"tl" (:@<.5ex>"tr"^-{\phi}:@<.5ex>"mr"^{Tf}:"b"^{y},:"ml"_{\Tensor\limits_ix_i}:@<.5ex>"b"_{f})} \]
commutes.

Now we suppose $(V,\tensor,T,\phi,\eta,\mu)$ is a symmetric monoidal monad and that $V$ is symmetric monoidal closed. Then given $X \in V$ and $(Y,y) \in V^T$, $[X,Y]$ obtains the structure of $T$-algebra, with the action corresponding under adjunction to the composite
\[ \xygraph{!{0;(2.5,0):} {T[X,Y]{\tensor}X}="p1" [r] {T[X,Y]{\tensor}TX}="p2" [r(.9)] {T([X,Y]{\tensor}X)}="p3" [r(.75)] {TY}="p4" [r(.5)] {Y}="p5"
"p1":"p2"^-{\id{\tensor}\eta}:"p3"^-{\phi}:"p4"^-{T\rev}:"p5"^-{y}}  \]
and we call this the \emph{pointwise $T$-algebra structure} on $[X,Y]$. The algebra axioms are easily verified. Moreover this assignment is functorial, that is to say, it is the object map of a functor
\[ [-,-] : \op V \times V^T \rightarrow V^T \]
which is compatible in the obvious way with the forgetful $U^T:V^T{\rightarrow}V$ and the original hom down in $V$. For all objects $X,Y \in V$ one may also define
\[ T_{X,Y} : [X,Y] \rightarrow [TX,TY] \]
corresponding under adjunction to the composite
\[ \xygraph{!{0;(2.5,0):} {[X,Y]{\tensor}TX}="p1" [r] {T[X,Y]{\tensor}TX}="p2" [r] {T([X,Y]{\tensor}X)}="p3" [r] {TY,}="p4"
"p1":"p2"^-{\eta{\tensor}\id}:"p3"^-{\phi}:"p4"^-{T\rev}} \]
and these maps give $T$ the structure of a $V$-functor. Moreover given $(X,x)$ and $(Y,y)$ in $V^T$, the composite
\[ \xygraph{!{0;(2,0):} {[X,Y]}="l" [r] {[TX,TY]}="m" [r] {[TX,Y]}="r" "l":"m"^-{T_{X,Y}}:"r"^-{[\id,y]}} \]
may be verified as the underlying map in $V$ of a $T$-algebra morphism, where $[X,Y]$ and $[TX,Y]$ are regarded as $T$-algebras via their pointwise structures.

Now let us suppose in addition that $V$ has equalisers. Then given $(X,x)$ and $(Y,y)$ in $V^T$, one can take the following equaliser
\begin{equation}\label{eq:alghom-def} \xygraph{{[(X,x),(Y,y)]}="eq" [r(2)] {[X,Y]}="l" [r(2)] {[TX,Y]}="r" [dl] {[TX,TY]}="m"
"eq":@{.>}"l"^-{e_{x,y}}(:"r"^-{[x,\id]},:"m"_{T_{X,Y}}:"r"_{[\id,y]})} \end{equation}
in $V$. By the recollections of the previous paragraph, this equaliser may also be regarded as living in $V^T$, and so the object $[(X,x),(Y,y)]$ as defined in equation(\ref{eq:alghom-def}) has a canonical $T$-algebra structure. The following result appears as theorem(2.2) of \cite{KockMM3}. We give a proof of it to illustrate the use of multicategories in the theory.
\begin{theorem}\cite{KockMM3}\label{thm:kock1}
Let $(V,\tensor,T,\phi,\eta,\mu)$ be a symmetric monoidal monad, $(V,\tensor)$ be symmetric monoidal closed with unit denoted as $I$, and let $V$ have equalisers. Then equation(\ref{eq:alghom-def}) defines the internal hom and $(TI,\mu_I)$ the unit of a closed structure on $V^T$.
\end{theorem}
\begin{proof}
In this proof we shall regard $V$ simultaneously as a representable multicategory and as a monoidal category without further comment. We shall denote by $V^T$ the symmetric multicategory of $T$-algebras. By proposition(\ref{prop:closedmcat->closedcat}) it suffices to exhibit $V^T$ as a closed multicategory with hom as given in equation(\ref{eq:alghom-def}), and to exhibit a universal $T$-algebra multimap $(){\rightarrow}(TI,\mu_I)$. But recall that for any $T$-algebra $(X,x)$, any multimap $(){\rightarrow}X$ down in $V$ satisfies the condition of a multimap of $T$-algebras. Moreover one has a universal multimap $u:(){\rightarrow}I$ down in $V$. From these observations together with the universal property of $(TI,\mu_I)$ as the free $T$-algebra on $I$, it follows easily that $Tu:(){\rightarrow}TI$ is universal in $V^T$.

Let $(Z,z)$ be in $V^T$. For any map $f:Z{\rightarrow}[X,Y]$ in $V$ let us write $\tilde{f}:(Z,X){\rightarrow}Y$ for the multimap in $V$ which corresponds to $f$ by closedness. It is easy to verify that $f$ is a $T$-algebra map, where the $T$-algebra structure on $[X,Y]$ is the pointwise one, iff $\tilde{f}$ is a $T$-algebra map in the first variable. It is also easy to verify that $[x,\id]f=[\id,y]T_{X,Y}f$ iff $\tilde{f}$ is a $T$-algebra morphism in the second variable. Thus by the universal property of the equaliser it follows that $\tilde{e}_{x,y}$ satisfies the universal property required of a right evaluation map in $V^T$.
\end{proof}
We impose one final additional assumption to obtain the induced tensor product of $T$-algebras. So we assume $(V,\tensor,T,\phi,\eta,\mu)$ satisfies the hypotheses of theorem(\ref{thm:kock1}), and in addition that $V^T$ has coequalisers. Then given a sequence $((X_i,x_i))_i$ of $T$-algebras, we may take the coequaliser
\begin{equation}\label{eq:algtensor-def}
\xygraph{{T\Tensor\limits_iTX_i}="l" [dr] {T^2\Tensor\limits_iX_i}="m" [ur] {T\Tensor\limits_iX_i}="r" [r(2)] {\Tensor\limits_i(X_i,x_i)}="c"
"l"(:@<1ex>"r"^-{T\Tensor\limits_ix_i},:"m"_{T\phi}:"r"_{\mu}:@<1ex>@{.>}"c"^-{q_{x_i}})} \end{equation}
in $V^T$, because the solid arrows here are clearly $T$-algebra morphisms between free $T$-algebras. The map $T\Tensor\limits_i\eta_{X_i}$ exhibits (\ref{eq:algtensor-def}) as a reflexive coequaliser. Note that while equation(\ref{eq:alghom-def}) may be regarded as living in both $V$ and $V^T$ since $U^T$ creates limits, equation(\ref{eq:algtensor-def}) is a coequaliser only in $V^T$: there is no reason in general for this coequaliser to be preserved by $U^T$. However in the motivating examples of this theory which involved finitary monads on $\Set$, this happened to be true because all finitary endofunctors of $\Set$ happen to preserve reflexive coequalisers, and thus the $U^T$ did create reflexive coequalisers in these cases.
\begin{theorem}\label{thm:smcc-from-monad}
Let $(V,\tensor,T,\phi,\eta,\mu)$ be a symmetric monoidal monad, $(V,\tensor)$ be symmetric monoidal closed with unit denoted as $I$, let $V$ have equalisers and $V^T$ have coequalisers. Then equation(\ref{eq:algtensor-def}) defines the tensor product and $(TI,\mu_I)$ the unit of a symmetric monoidal closed structure on $V^T$.
\end{theorem}
\begin{proof}
We continue to regard $V$ simultaneously as a representable multicategory and as a monoidal category. By theorem(\ref{thm:kock1}), lemma(\ref{lem:closed->representable}) and proposition(\ref{prop:smultcat->smoncat}), it suffices to exhibit universal multimaps
\[ \begin{array}{c} {((X_i,x_i))_i \rightarrow \Tensor\limits_i(X_i,x_i).} \end{array} \]
Let $(Z,z)$ be a $T$-algebra. For any $T$-algebra map
\[ \begin{array}{c} {f:(T\Tensor\limits_iX_i,\mu){\rightarrow}(Z,z),} \end{array} \]
denote by $\tilde{f}:(X_i)_i{\rightarrow}Z$ the multimap down in $V$ which corresponds to $f$. Then one may easily verify that
\[ \begin{array}{c} {fT(\Tensor\limits_ix_i)=f{\mu}T\phi} \end{array} \]
iff $\tilde{f}$ satisfies the conditions of a multi-$T$-algebra morphism. Thus $\tilde{q}_{x_i}$ is a universal multimap as required. 
\end{proof}
\begin{remark}\label{rem:tensor-of-frees}
In the particular case where the $((X_i,x_i))_i$ are free, say $(X_i,x_i)=(TZ_i,\mu_{Z_i})$, then one may take $q_{x_i}$ in (\ref{eq:algtensor-def}) to be given by the composite
\[ \xygraph{!{0;(2,0):} {T\Tensor\limits_iTZ_i}="l" [r] {T^2\Tensor\limits_iZ_i}="m" [r] {T\Tensor\limits_iZ_i}="r" "l":@<1ex>"m"^-{T\phi}:@<1ex>"r"^-{\mu}} \]
because then the maps 
\[ \xygraph{!{0;(2.5,0):} {T\Tensor\limits_iT^2Z_i}="l" [r] {T\Tensor\limits_iTZ_i}="m" [r] {T\Tensor\limits_iZ_i}="r" "r":@<-1ex>"m"_-{T\Tensor\limits_i\eta_{Z_i}}:@<-1ex>"l"_-{T\Tensor\limits_iT\eta_{Z_i}}} \]
provide a splitting exhibiting our suggested $q_{x_i}$ as a split coequaliser. Thus one may regard $F^T$ as a strict monoidal functor. In the approach to monoidal monad theory of Brian Day, this is the starting point -- one defines the tensor product of free algebras in this way, and then extends the definition of the induced tensor product to all of $V^T$ using convolution and the density of $V_T$ within $V^T$. See \cite{DayDPres} for more details.
\end{remark}
%

\section{Symmetric multicategories of enriched graphs}\label{sec:multicategories-of-graphs}

\subsection{Overview}
In this section we exhibit $\ca GV$ as category of linear maps underlying a symmetric multicategory, by applying a 2-functor
\[ \ca F : \CAT/\Set \rightarrow \SYMMULT \]
to $\ca GV$, regarded as over $\Set$ in the usual way by the functor which sends an enriched graph to its set of objects. While an object of $\CAT/\Set$ is really a functor $A \rightarrow \Set$, we shall throughout this section suppress mention of the functor, denote by $a_0$ the underlying set of $a \in A$, and refer to an element $x$ of $a_0$ as an element or object of $a$, sometimes writing $x \in a$ to denote this.

\subsection{The object map of $\ca F$}\label{ssec:obj-F}
If $X$ is a symmetric multicategory and $(x_1,...,x_n)$ a sequence of objects of $X$, then taking homs produces a functor
\[ X((x_i)_i,-) : \lin(X) \rightarrow \Set. \]
The cartesian monoidal structure of $\Set$ enables us to view it as a symmetric multicategory by applying $U$, and we denote this symmetric multicategory by $\SetAsMCat$. Obviously $\lin(\SetAsMCat)=\Set$. Thus it makes sense to ask whether the above hom functor for $X$ may be viewed as a symmetric multifunctor into $\SetAsMCat$. In order to do so, for each $f:(a_1,...,a_m){\rightarrow}b$  in $X$, it is necessary to define a function
\[ X((x_i)_i,a_1) \times ... \times X((x_i)_i,a_m) \rightarrow X((x_i)_i,b) \]
using composition in $X$. So one is tempted to guess that this function is defined by
\[ (g_1,....,g_m) \mapsto f(g_1,...,g_m). \]
While $f(g_i)_i$ is a well-defined multimap, it doesn't live in the correct hom-set, but instead is of the form
\[ (x_1,...,x_n,......,x_1,...,x_n) \rightarrow b \]
where the domain sequence is the concatenation of $m$-copies of $(x_1,..,x_n)$. There are various ways one can imagine to get around this ``problem''. However there is one case where there is no problem, that is when $(x_1,...,x_n)$ happens to be the empty sequence $()$. Thus homming out of the empty sequence produces both a functor and a symmetric multifunctor
\[ \begin{array}{lccr} {X((),-) : \lin(X) \rightarrow \Set} &&& {X((),-) : X \rightarrow \SetAsMCat{,}} \end{array} \]
and one obtains the former from the latter by applying the 2-functor $\lin$.

In this way given a symmetric multicategory $X$ and an object $x$ therein, an element $z$ of $x$ \emph{is} a multimap $z:(){\rightarrow}x$ in $X$. Given any multimap \[ f:(a_1,...,a_n){\rightarrow}b \] in $X$, the assignment
\[ \begin{array}{lcr} {z \in (a_{i,0})_i} & \mapsto & {f(z_1,...,z_n).} \end{array} \]
describes the multimap $X((),f)$ of $\SetAsMCat$. Given a choice $1{\leq}i^*{\leq}n$ of input variable, and a choice $z \in (a_{i,0})_{i{\neq}i^*}$ of element for all but the $i^*$-th variable, one can define a linear map
\[ f_z : a_{i^*} \rightarrow b \]
as the composite
\[ f(z_1,...,z_{i^*{-}1},1_{a_{i^*}},z_{i^*{+}1},...,z_n) \]
that is, as the result of fixing all but the $i^*$-th variable at the values $z$ and applying $f$.
\begin{definition}
Let $f:(a_1,...,a_n){\rightarrow}b$ be a multimap in a symmetric multicategory $X$. The linear maps $f_z$ obtained from $f$ in the manner just described are called the \emph{linear parts} of $f$. The multimap $X((),f)$ of $\SetAsMCat$ is called the \emph{underlying multifunction of $f$} and is denoted as $f_0$.
\end{definition}
\noindent Note in particular that if $f$ is itself linear (ie when $n=1$), then $i^*$ can only be $1$ and so $z$ can only be the empty sequence $()$ and so the only linear part of $f$ is $f$ itself, that is, $f=f_{()}$. The way in which the linear parts of multimaps in $X$ get along with the composition and symmetric group actions of $X$ is expressed by the following lemma, whose proof is a simple exercise in the definitions just given.
\begin{lemma}\label{lem:alg-linparts}
Given multimaps
\[ \begin{array}{lcr} {f:(a_1,...,a_n){\rightarrow}b} && {g:(c_{i1},...,c_{im_i}){\rightarrow}a_i} \end{array} \]
in a symmetric multicategory $X$, choices $1{\leq}i^*{\leq}n$ and $1{\leq}j^*{\leq}m_{i^*}$ of input variables, and sequences of elements
\[ \begin{array}{lcr} {z \in (a_{i,0})_{i{\neq}i^*}} && {w \in (c_{ij,0})_{(i,j){\neq}(i^*,j^*)}} \end{array} \]
and a permutation $\sigma \in \Sigma_n$.
Then the formulae
\[ \begin{array}{lcr} {(f\sigma)_z=f_{z\sigma^{-1}}} && {(f(g_i)_i)_w = f_{g(w_{ij})_{i{\neq}i^*}}g_{(w_{(i^*j)})_{j{\neq}j^*}}.} \end{array} \]
are valid.
\end{lemma}
\noindent To summarise, for any multicategory $X$ and any multimap $f$ therein, one has an underlying multimap $f_0$ in $\SetAsMCat$ and linear parts the $f_z$ just defined, and one can be quite precise about how the algebra of $X$ is reflected in this data: $f \mapsto f_0$ is the multimap assignment of a symmetric multifunctor, namely $X((),-)$, and lemma(\ref{lem:alg-linparts}) describes what happens at the level of the linear parts.

In the case $X = \SetAsMCat$, $f$ and $f_0$ are the same, and moreover when $n{\neq}0$ the linear parts of $f$ also determine $f$ uniquely. When $n=0$ there are no linear parts but such an $f$ is an element of $b$. By contrast it can happen that the underlying multifunction and linear parts together tell us very little about a given multimap in a given multicategory. Consider for instance the case where $X$ is the symmetric multicategory whose object set is\[ \{0,...,n{+}1\} \] where $n>0$, and the hom $X((0,...,n),n{+}1)$ is some arbitrary set $Z$. Suppose that $X$ is freely generated as a symmetric multicategory with this data. One way to describe the resulting $X$, by the yoneda lemma, is that a morphism $X{\rightarrow}Y$ in $\SYMMULT$ amounts to a choice \[ a_0,...,a_n,b \] of objects of $Y$ together with a function $Z{\rightarrow}Y((a_i)_i,b)$. Note that in $X$ there are no multimaps out of the empty sequence, thus the set of elements of any $x \in X$ is empty. Thus the underlying multifunction of any $f \in Z$ is the empty one and $f$ has no linear parts.

Given a category $A$ over $\Set$, the symmetric multicategory $\ca FA$ is defined so that one recovers the original $A{\rightarrow}\Set$ as the $\Set$-valued hom functor $\ca FA((),-)$, and the multimaps of $\ca FA$ are determined uniquely by their underlying multifunctions and linear parts. The explicit definition is
\begin{definition}
An object of $\ca FA$ is an object of $A$. A multimap \[ f : (a_1,...,a_n) \rightarrow b \] in $\ca FA$ is given by the data
\[ \begin{array}{lccr} {f_0 : (a_{1,0},...,a_{n,0}) \rightarrow b_0} &&& {f_z : a_{i^*} \rightarrow b} \end{array} \]
where for $z \in (a_{i,0})_{i{\neq}i^*}$, $f_z$ is a morphism of $A$. This data must satisfy $(f_z)_0=(f_0)_z$ for all $z$.
\end{definition}
\noindent Note that when $n=0$ such a multimap is just an object of $b$, and when $n=1$, $f$ may be identified with the morphism $f_{()}$ of $A$. With this understood, we define the identities of $\ca FA$ to be those of $A$.

Let $f:(a_1,...,a_n) \rightarrow b$ be a morphism of $\ca FA$ and $\sigma \in \Sigma_n$. Then we define
\[ f\sigma : (a_{\sigma{i}})_i \rightarrow b \]
by
\[ \begin{array}{lccr} {(f\sigma)_0=f_0\sigma} &&& {(f\sigma)_z=f_{z\sigma^{-1}}} \end{array} \]
for all $1{\leq}i^*{\leq}n$ and $z \in (a_{\sigma{i},0})_{i{\neq}i^*}$. By definition and lemma(\ref{lem:alg-linparts}) we have the calculation
\[ ((f\sigma)_z)_0 = (f_{z\sigma^{-1}})_0 = (f_0)_{z\sigma^{-1}} = (f_0\sigma)_z = ((f\sigma)_0)_z \]
exhibiting the well-definedness of $f\sigma$.

Let $f$ be as above and $g_i:(c_{i1},...,c_{im_i}){\rightarrow}a_i$ be morphisms of $\ca FA$ for $1{\leq}i{\leq}n$. Then we define the composite
\[ f(g_i)_i : (c_{ij})_{ij} \rightarrow b \]
in $\ca FA$ by
\[ \begin{array}{lccr} {(f(g_i)_i)_0 = f_0(g_{i,0})_i} &&& {(f(g_i)_i)_w = f_{g(w_{ij})_{i{\neq}i^*}}g_{(w_{(i^*j)})_{j{\neq}j^*}}} \end{array} \]
for all choices $1{\leq}i^*{\leq}n$ and $1{\leq}j^*{\leq}m_{i^*}$ of input variables and sequences of elements $w \in (c_{ij,0})_{(i,j){\neq}(i^*,j^*)}$. By a similar calculation as with $f\sigma$ using the definitions and lemma(\ref{lem:alg-linparts}), one may verify that
\[ ((f(g_i)_i)_w)_0 = ((f(g_i)_i)_0)_w \]
and so exhibit $f(g_i)_i$ as well-defined.
\begin{proposition}
Let $A$ be a category over $\Set$. With the data, symmetric group actions and compositions just defined, $\ca FA$ forms a symmetric multicategory. The $\Set$-valued hom $\ca FA((),-)$ is the original $A{\rightarrow}\Set$.
\end{proposition}
\begin{proof}
The symmetric group actions $\ca FA$ are functorial because they are so at the level of the underlying multifunctions. The unit and associative laws of composition follow from those of $\SetAsMCat$ and $A$, and lemma(\ref{lem:alg-linparts}). Equivariance of composition in $\ca FA$ follows from that of $\SetAsMCat$. As explained already, one identifies any linear map $f$ in $\ca FA$ with the morphism $f_{()}$ of $A$, and under this identification composition of linear maps in $\ca FA$ corresponds to composition in $A$. Moreover under this identification of linear maps, for any multimap $f$ in $\ca FA$, one may identify the $f_z$'s of its definition as the corresponding linear parts of $f$. This last is just the identification of the $\Set$-valued hom $\ca FA((),-)$ with the original $A{\rightarrow}\Set$.
\end{proof}
\begin{example}\label{ex:FGV}
Let $V$ be a category and regard the category $\ca GV$ as over $\Set$ in the usual way. For a multimap in $\ca F(\ca GV)$
\[ f : (A_1,...,A_n) \rightarrow B \]
of $V$-graphs, the underlying multifunction $f_0$ determines the object maps of the linear parts $f_z$ where $z \in (A_i)_{i{\neq}i^*}$. Thus in addition to $f_0$, the data for $f$ involves hom maps
\[ f_{z,a,b} : A_{i^*}(a,b) \rightarrow B(f(z|_{i^*}a),f(z|_{i^*}b)) \]
in $V$ for all $1{\leq}i^*{\leq}n$, $z \in (A_i)_{i{\neq}i^*}$ and $a,b \in A_{i^*}$, and this data satisfies no further conditions. The case $V=\Set$ and $n=2$ is instructive, for then the hom maps involve assignments of the form
\[ \begin{array}{lcr} {(a_1,\alpha_2:a_2{\rightarrow}a_2')} & \mapsto & {f(a_1,\alpha_2):f(a_1,a_2){\rightarrow}f(a_1,a_2')} \\
{(\alpha_1:a_1{\rightarrow}a_1',a_2)} & \mapsto & {f(\alpha_1,a_2):f(a_1,a_2){\rightarrow}f(a_1',a_2)} \end{array} \]
where $a_1,a_1',\alpha_1$ are vertices and edges from $A_1$, and $a_2,a_2',\alpha_2$ are vertices and edges from $A_2$. Thus from the edges $\alpha_1$ and $\alpha_2$, one obtains the following square
\[ \xygraph{!{0;(2,0):(0,.5)::} {f(a_1,a_2)}="tl" [r] {f(a_1',a_2)}="tr" [d] {f(a_1',a_2')}="br" [l] {f(a_1,a_2')}="bl"
"tl"(:"tr"^-{f(\alpha_1,a_2)}:"br"^{f(a_1',\alpha_2)},:"bl"_{f(a_1,\alpha_2)}:"br"_-{f(\alpha_1,a_2')})} \]
in the graph $B$. In a similar way for general $n$, $f$ produces an $n$-dimensional hypercube in the graph $B$ from a given $n$-tuple of edges $(\alpha_i)_i$ from the $A_i$.
\end{example}
\begin{example}
Regard $\Cat$ as over $\Set$ via the functor which sends a category to its set of objects. Then a multimap of categories
\[ f : (A_1,...,A_n) \rightarrow B \]
in $\ca F(\Cat)$ amounts to a multimap in $\ca F(\ca G\Set)$ as in example(\ref{ex:FGV}), together with the requirement that the hom functions be functorial in each variable separately. Note that the hypercubes in the category $B$ described in example(\ref{ex:FGV}) will not necessarily commute in general.
\end{example}
\begin{example}
Let $X$ be a set and $M$ the monoid of endofunctions of $X$. The monoid $M$ acts on $X$ by evaluation, and regarding $M$ as a one object category this action may be expressed as a functor $M{\rightarrow}\Set$ in which the unique object of $M$ is sent to the set $X$. The symmetric multicategory $\ca FM$ has one object, thus it is just an operad of sets. One may easily verify directly that $\ca FM$ is in fact the endomorphism operad of $X$.
\end{example}
\begin{example}
For $R$ a commutative ring and regard $\RMod$ as over $\Set$ via the forgetful functor. Then $\ca F(\RMod)$ is the symmetric multicategory of $R$-modules and $R$-multilinear maps.
\end{example}
%

\subsection{2-functoriality of $\ca F$}\label{ssec:2functoriality-F}
We shall define $\ca F$ as the composite 2-functor
\[ \xygraph{!{0;(3,0):} {\CAT/\Set}="l" [r] {\SYMMULT/\SetAsMCat}="m" [r] {\SYMMULT}="r" "l":"m"^-{\ca F_1}:"r"^-{\dom}} \]
where $\dom$ 
is obtained by talking the domain of a map into $\SetAsMCat$, and the 2-functor $\ca F_1$ is provided by
\begin{lemma}
The assignment
\[ \begin{array}{lcr} A{\rightarrow}\Set & \mapsto & {\ca FA((),-):\ca FA{\rightarrow}\SetAsMCat} \end{array} \]
describes the object map of a 2-functor which we denote as $\ca F_1$.
\end{lemma}
\begin{proof}
Let $F:A{\rightarrow}B$ be a functor over $\Set$. We define $\ca F_1(F)$ as follows. Its object map is that of $F$. Given $f : (a_1,...,a_n){\rightarrow} b$ in $\ca FA$, we define
\[ \begin{array}{lccr} {\ca F_1(F)(f)_0 = f_0} &&& {\ca F_1(F)(f)_z = F(f_z)} \end{array} \]
Let $\phi:F{\rightarrow}G$ be a natural transformation over $\Set$, then in view of the identification of the linear maps of $\ca FB$ with the morphisms of $B$, we may define the components of $\ca F_1(\phi)$ to be those of $\phi$. It is trivial to verify that this one and 2-cell map for $\ca F_1$ is well-defined and 2-functorial.
\end{proof}
%

\subsection{Symmetric monoidal monads from monads over $\Set$}\label{ssec:monads-from-F}
Under appropriate conditions, a monad over $\Set$ has a canonical symmetric monoidal structure because the 2-functor $\ca F$ gets along well with monad theory. Before explaining this we must first clarify what Eilenberg-Moore objects are in the slices of a 2-category.
\begin{lemma}
Let $\ca K$ be a 2-category, $X \in \ca K$ and
\[ \dom : \ca K/X \rightarrow \ca K \]
be the 2-functor which on objects takes the domain of a morphism into $X$. Then $\dom$ creates any Eilenberg-Moore objects that exist in $\ca K$.
\end{lemma}
\begin{proof}
Let $f:A{\rightarrow}X$ and $(T,\eta,\mu)$ be a monad on $f$ in $\ca K/X$, or in other words a monad on $A$ in $\ca K$ which satisfies $fT=f$ and $f\eta=\id_f=f\mu$. Then any algebra of $T$ in $\ca K$, that is to say $z:Z{\rightarrow}A$ and $\alpha:Tz{\rightarrow}z$ satisfying the usual axioms, automatically lives in $\ca K/X$: the underlying object is $fz$ and $f\alpha=\id$ by the unit law of $(z,\alpha)$ since $f\eta=\id$. Thus the Eilenberg-Moore object $(U^T,\tau)$ of $T$ in $\ca K$, when it exists, may be regarded as living in $\ca K/X$ with underlying object $fU^T$. By the one-dimensional part of the universal property of $(U^T,\tau)$ in $\ca K$, for each $(z,\alpha)$ one has a unique $z':Z{\rightarrow}A^T$ such that $U^Tz'=z$, but then post-composing this with $f$ shows that $z'$ also lives over $X$, and so $(fU^T,\tau)$ enjoys the one-dimensional part of the universal property of an Eilenberg-Moore object for $T$ in $\ca K/X$. Given a morphism $\phi:(z,\alpha) \rightarrow (z',\alpha')$ of $T$-algebras in $\ca K/X$, the two-dimensional part of the universal property of $(U^T,\tau)$ in $\ca K$ gives us a unique $\phi':z{\rightarrow}z'$ such that $U^T\phi'=\phi$, so $fU^T\phi'=f\phi=\id$ since $\phi$ lived over $X$ in the first place, and so $(fU^T,\tau)$ enjoys the one-dimensional part of the required universal property.
\end{proof}
\noindent So there is in particular no real distinction between Eilenberg-Moore objects in $\CAT/\Set$ and those in $\CAT$. The sense in which $\ca F$ gets along well with monad theory is described by
\begin{lemma}\label{lem:F-EM}
$\ca F$ preserves Eilenberg-Moore objects.
\end{lemma}
\begin{proof}
Let $X$ be a category over $\Set$ and $T$ a monad on $A$ over $\Set$. We must compare the symmetric multicategories $\ca F(X^T)$ and $\ca F(X)^{\ca F(T)}$, this latter being formed by applying $\ca F$ to the monad and taking Eilenberg-Moore objects in $\SYMMULT$, which we know how to do by proposition(\ref{prop:EM-SYMMULT}). By definition each of these multicategories has the same objects, namely algebras of the monad $T$. A multimap
\[ ((x_1,a_1),...,(x_n,a_n)) \rightarrow (y,b) \]
of $\ca F(A^T)$ consists of a multimap $f:(x_{i,0})_i \rightarrow y_0$ in $\SetAsMCat$, together with $f_z:x_{i^*} \rightarrow y$ in $X^T$ for all $1{\leq}i^*{\leq}n$ and $z \in (x_{i,0})_{i{\neq}i^*}$. This is clearly the same as a multimap $f':(x_i)_i \rightarrow y$ in $\ca F(X)$ whose linear parts are $T$-algebra morphisms. On the other hand a multimap $((x_i,a_i))_i \rightarrow (y,b)$ in $\ca F(X)^{\ca F(T)}$ is by definition a multimap $f'$ such that
\[ f'(\id_{x_i}|_{i^*}a_{i^*})_i = bT(f')(\eta_{x_i}|_{i^*}\id_{Tx_{i^*}})_i \]
as multimaps of $\ca F(X)$. But two multimaps in $\ca F(X)$ are equal iff they have the same underlying multifunction and the same linear parts. The underlying multifunctions of the multimaps on both sides of the previous equation is just the identity since $T$ is a monad over $\Set$ and so the underlying function of a given $T$-action must be the identity by the unit law. The equality of linear parts expressed by the equation says exactly that these linear parts are $T$-algebra maps. The identifications of objects and multimaps between $\ca F(X^T)$ and $\ca F(X)^{\ca F(T)}$ that we have just made clearly commute in the obvious way with the forgetful $U^T:X^T{\rightarrow}X$, that is to say, we have shown that the obstruction map $\ca F(X^T) \rightarrow \ca F(X)^{\ca F(T)}$ is an isomorphism as required.
\end{proof}
With these details in hand we come to the main result of this section, which explains how monads over $\Set$ naturally give rise to symmetric monoidal monads, and when this happens, the two obvious ways of regarding the algebras as multicategories coincide.
\begin{theorem}\label{thm:MndOverSet->SMultCat}
Let $A$ be a category over $\Set$ and $T$ a monad on $A$ over $\Set$. Suppose that $\ca FA$ is a representable multicategory.
\begin{enumerate}
\item  $T$ is canonically a symmetric monoidal monad relative to the induced symmetric monoidal structure on $A$.
\item  One has an isomorphism $\ca F(A^T) \iso UA^{UT}$ of symmetric multicategories. 
\end{enumerate}
\end{theorem}
\begin{proof}
Since $U$ is 2-fully-faithful by proposition(\ref{prop:smultcat->smoncat}) and $T$ is just $\ca FT$ restricted to the linear maps in $\ca FA$, one may regard $T$ in a unique way as a symmetric monoidal monad such that $UT=\ca FT$. The second part follows since $\ca F$ preserves Eilenberg-Moore objects by lemma(\ref{lem:F-EM}).
\end{proof}
As we shall see in the next section, it is easy to isolate conditions on $V$ so that the symmetric multicategory $\ca F(\ca GV)$ is closed and representable. For such $V$ and given a monad $T$ on $\ca GV$ over $\Set$, theorem(\ref{thm:MndOverSet->SMultCat}) tells us that the multicategory $\ca F(\ca G(V)^T)$ arises from a symmetric monoidal structure on the monad $T$. Thus one may use theorem(\ref{thm:smcc-from-monad}) to study the representability of $\ca F(\ca G(V)^T)$, and so ultimately produce a canonical symmetric monoidal closed structure on $\ca G(V)^T$ in a very general way. Such technique will apply in particular to categories of algebras of higher operads.

\section{The free tensor products}\label{sec:free}

\subsection{Overview}
Having assembled together the necessary technology in sections(\ref{sec:monoidal-monads}) and (\ref{sec:multicategories-of-graphs}) we are now in a position to describe our generalisations of the free product of categories recalled in section(\ref{sec:funny}), and exhibit some of their basic properties.

\subsection{Free products of enriched graphs}\label{ssec:fp-GV}
We begin by considering the multicategory $\ca F(\ca GV)$ for some fixed category $V$, and we shall use the notation introduced in example(\ref{ex:FGV}) when specifying multimaps therein. Supposing first that $V$ has small products, one can for each pair $A$, $B$ of $V$-enriched graphs define $[A,B]$ in $\ca GV$ to have as objects, morphisms $f:A{\rightarrow}B$ of enriched graphs, and homs given by
\begin{equation}\label{eq:GV-hom} \begin{array}{c} {[A,B](f,g) = \prod\limits_{a{\in}A} B(fa,ga).} \end{array} \end{equation}
The right evaluation $\rev_{A,B}:([A,B],A){\rightarrow}B$ is specified as follows
\[ \begin{array}{lcccr} {\rev_0(f,a) = fa} && {\rev_{f,a_1,a_2} = f_{a_1,a_2}} && {\rev_{a,f,g} = p_a} \end{array} \]
where $f,g$ are $V$-graph morphisms $A{\rightarrow}B$, $a,a_1,a_2 \in A$ and $p_a$ is the $a$-th projection of the product from (\ref{eq:GV-hom}). Recall that $f_{a_1,a_2}$ is our notation for the corresponding hom map of the $V$-graph morphism $f$.
\begin{proposition}\label{prop:closed-FGV}
Let $V$ be a category with small products. With the hom and right evaluation multimaps just given, $\ca F(\ca GV)$ is a closed multicategory.
\end{proposition}
\begin{proof}
Given $F$ in $\ca F(\ca GV)$ as in
\[ \begin{array}{lccr} {F:(C_1,...,C_n,A) \rightarrow B} &&& {G:(C_1,...,C_n) \rightarrow [A,B]} \end{array} \]
we must exhibit a unique $G$ as above so that $F=\rev(G,1_A)$. Let us describe the object map of $G$. Given $x \in (C_i)_i$ one must define a morphism
\[ G(x_i)_i : A \rightarrow B \]
of $V$-graphs, and the object and hom maps are defined by
\[ \begin{array}{lccr} {G(x_i)_i(a) = F(x|_{n{+}1}a)} &&& {(G(x_i)_i)_{a_1,a_2} = F_{x,a_1,a_2}.} \end{array} \]
As for the hom maps of $G$ for $x \in (C_i)_{i{\neq}i^*}$ and $y,z \in C_{i^*}$, we define $G_{x,y,z}$ as the unique map in $V$ satisfying
\[ p_aG_{x,y,z} = F_{x|_{n{+}1}a,y,z} \]
for all $a \in A$ where $p_a$ is the projection from the product (\ref{eq:GV-hom}). It is immediate from the definitions just given that $F=\rev(G,1)$. Conversely, the closedness of $\SetAsMCat$ and the definition of the object map of $\rev$ ensures that $F=\rev(G,1)$ determines the object maps of the $G(x_i)_i$. For $x \in (C_i)_i$ and $a_1,a_2 \in A$, observing the corresponding homs on both sides of $F=\rev(G,1)$ one sees that the hom maps of the $G(x_i)_i$ are also determined uniquely by that formula. Finally given $x \in (C_i)_{i{\neq}i^*}$, $a \in A$, and $y,z \in C_{i^*}$, observing the corresponding homs on both sides of $F=\rev(G,1)$ one sees that the morphisms $p_aG_{x,y,z}$ are also determined by $F=\rev(G,1)$, and so we have indeed exhibited $G$ as the unique solution of $F=\rev(G,1)$ as required.
\end{proof}
Suppose now that in addition to having small products, that $V$ also has finite coproducts. Then given $V$-graphs $(A_1,...,A_n)$ one can define a $V$-graph \[ \begin{array}{c} {\FreeProd\limits_iA_i} \end{array} \] with object set the cartesian product of the $A_{i,0}$, and homs given by
\begin{equation}\label{GV-tensor} \begin{array}{c} {\left(\FreeProd\limits_iA_i\right)(a,b) = \left\{ \begin{array}{lll} {\coprod_i A_i(a_i,b_i)} && {\textnormal{if $a_i{=}b_i$ for all $i$.}} \\ {A_j(a_j,b_j)} && {\textnormal{if $a_i{=}b_i$ for all $i$ except $i{=}j$.}} \\ {\emptyset} && {\textnormal{otherwise.}} \end{array} \right.} \end{array} \end{equation}
and we define a multimap
\[ \begin{array}{c} {\alpha_{A_i} : (A_i)_i \rightarrow \FreeProd\limits_iA_i} \end{array} \]
whose object map is the identity, and whose hom map
\[ \begin{array}{c} {A_{i^*}(y,z) \rightarrow (\FreeProd\limits_i A_i)(x|_{i^*}y,x|_{i^*}z)} \end{array} \]
corresponding to $x \in (A_i)_{i{\neq}i^*}$, $y,z \in A_{i^*}$ is the identity if $y \neq z$, and the $i^*$-th coproduct inclusion otherwise.
\begin{proposition}\label{prop:FGV-rep}
Let $V$ be a category with small products and finite coproducts. Then $\ca F(\ca GV)$ is representable and the multimaps $\alpha_{A_i}$ just defined are universal.
\end{proposition}
\begin{proof}
By lemma(\ref{lem:closed->representable}) and proposition(\ref{prop:closed-FGV}) it suffices to show that the $\alpha_{A_i}$ are universal. Given $F$ in $\ca F(\ca GV)$ as in
\[ \begin{array}{lccr} {F : (A_i)_i \rightarrow B} &&& {G : \FreeProd\limits_i A_i \rightarrow B} \end{array} \]
we must exhibit $G$ as above unique so that $G\alpha=F$. Clearly this equation forces the object map of $G$ to be that of $F$, and the hom map $G_{a,b}$ to be: (1) the unique map such that $G_{a,b}c_i = F_{a{\neg}i,a_i,b_i}$ when $a{=}b$ where $a{\neg}i$ is the sequence $a$ with $i$-th coordinate removed, (2) $F_{a{\neg}j,a_j,b_j}$ when $a_i{=}b_i$ for all $i$ except $i{=}j$, and (3) the unique map out of $\emptyset$ otherwise.
\end{proof}
\begin{corollary}
If $V$ is a category with small products and finite coproducts, then $\FreeProd$ and the hom described in equation(\ref{eq:GV-hom}) give $\ca GV$ a symmetric monoidal closed structure, and for this structure $U\ca GV = \ca F(\ca GV)$.
\end{corollary}
\begin{definition}
If $V$ is a category with small products and finite coproducts, then the tensor product $\FreeProd$ is called the \emph{free product} of enriched graphs.
\end{definition}
\noindent Note that the unit for the free product, that is to say the nullary case, is just the $V$-graph $0$ which represents the forgetful $\ca GV{\rightarrow}\Set$: it has one object and its unique hom is $\emptyset$.
\begin{example}\label{ex:free-product-Graph}
Unpacking the binary free product in the case $V=\Set$ gives the simpler graph-theoretic analogue of the funny tensor product of categories. As with categories, $A \BinFreeProd B$ and the cartesian product $A \times B$ have the same objects, but their edges are different. There are two types of edges of $A \BinFreeProd B$:
\[ \begin{array}{lccr} {(a,\beta) : (a,b_1) \rightarrow (a,b_2)} &&& {(\alpha,b) : (a_1,b) \rightarrow (a_2,b)} \end{array} \]
where $a$ and $\alpha:a_1{\rightarrow}a_2$ are in $A$, and $b$ and $\beta:b_1{\rightarrow}b_2$ are in $B$. In particular given such an $\alpha$ and $\beta$, the dotted arrows in
\[ \xygraph{!{0;(3,0):(0,.333)::} {(a_1,b_1)}="tl" [r] {(a_1,b_2)}="tr" [d] {(a_2,b_2)}="br" [l] {(a_2,b_1)}="bl" "tl"(:@{.>}"tr"^-{(a_1,\beta)}:@{.>}"br"^{(\alpha,b_2)},:@{.>}"bl"_{(\alpha,b_1)}:@{.>}"br"_-{(a_2,\beta)},:"br"|{(\alpha,\beta)})} \]
indicate some edges which one can build from them in $A \BinFreeProd B$, whereas the solid diagonal edge is what one has in $A \times B$.
\end{example}
There is an importance difference between the free product of graphs and that of categories. That is for graphs there is no sensible comparison map
\[ A \BinFreeProd B \rightarrow A \times B \]
as in the $\Cat$ case.

\subsection{Free products of monad algebras}\label{ssec:fp-monad-algebras}
Assembling together the free product of enriched graphs just exhibited, monoidal monad theory and theorem(\ref{thm:MndOverSet->SMultCat}), we can now exhibit analogous tensor products of algebras of monads on $\ca GV$ over $\Set$.
\begin{theorem}\label{thm:free-product-of-algebras}
If $T$ is an accessible monad on $\ca GV$ over $\Set$ and $V$ be locally presentable, then $\ca F(\ca G(V)^T)$ is closed and representable. 
\end{theorem}
\begin{proof}
$\ca GV$ is locally presentable by corollary(5.12) of \cite{EnHopII} or by \cite{KL-NiceVCat}, and since $T$ is accessible $\ca G(V)^T$ is locally presentable also and so $T$ satisfies the hypotheses of theorems(\ref{thm:kock1}) and (\ref{thm:smcc-from-monad}) from which the result follows.
\end{proof}
\begin{definition}
The induced tensor product on $\ca G(V)^T$ by theorem(\ref{thm:free-product-of-algebras}) is called the \emph{free product} on $\ca G(V)^T$ and is also denoted as $\FreeProd$.
\end{definition}
\noindent Equations(\ref{eq:alghom-def}) and (\ref{eq:algtensor-def}) of section(\ref{ssec:smonmnd}), together with the explicit description of free products and the associated internal hom for $\ca GV$, gives one explicit formulas for the tensors and homs provided by theorem(\ref{thm:free-product-of-algebras}). In particular one has the following simple observation.
\begin{remark}\label{rem:unit-of-free-product}
For $T$ as in theorem(\ref{thm:free-product-of-algebras}) the unit of the induced monoidal structure is $(T0,\mu_0)$, that is to say, the free $T$-algebra on the graph $0$ with one object and initial hom.
\end{remark}
\noindent Moreover unpacking the equations (\ref{eq:alghom-def}) and (\ref{eq:algtensor-def}) in the case where $V=\Set$ and $T$ is the category monad $\ca T_{{\leq}1}$, one recovers the tensor and hom of the funny tensor product of categories. Theorem(\ref{thm:free-product-of-algebras}) applies to any higher operad, and so any category of algebras of a higher operad has such a symmetric monoidal closed structure.

\subsection{Comparing the free and cartesian products}\label{ssec:free-vs-cartesian} From the above discussion at the end of section(\ref{ssec:fp-GV}) we see that the category monad $T=\ca T_{{\leq}1}$ on $\ca G\Set$ is ``better'' than the identity monad on $\ca G\Set$, because in $\ca G\Set^T=\Cat$ one has natural identity on objects comparison maps mediating between the free and cartesian product. We shall now isolate which formal property of $T$ gives rise to these comparison maps in general.

Suppose that $T$ is a monad over $\Set$ on $\ca GV$ with this situation satisfying the conditions of theorem(\ref{thm:free-product-of-algebras}) so that one has free and cartesian products in $\ca G(V)^T$. Suppose in addition that one has for all $A$ and $B$ in $\ca G(V)^T$ some map
\[ \kappa_{A,B} : A \BinFreeProd B \rightarrow A \times B. \]
We don't assume anything about the naturality of the $\kappa_{A,B}$ or even that they are identities on objects. Then in particular putting $A=1$ and $B=(T0,\mu_0)$, in view of remark(\ref{rem:unit-of-free-product}) one obtains by composition with the appropriate coherence isomorphisms, a $T$-algebra map
\[ e : 1 \rightarrow (T0,\mu_0). \]
Of course as a map out of a terminal object $e$ is a split monomorphism. In fact in this case it must be an isomorphism. To see that the composite
\[ \xygraph{{T0}="l" [r] {1}="m" [r] {T0}="r" "l":"m":"r"^-{e}} \]
is the identity, by the universal property of $T0$ as the free $T$-algebra on $0$, it suffices to show that the composite morphism
\[ \xygraph{{0}="p1" [r] {T0}="p2" [r] {1}="p3" [r] {T0}="p4" "p1":"p2"^-{\eta_0}:"p3":"p4"^-{e}} \]
in $\ca GV$ is $\eta_0$, and this follows since $0$ is the initial $V$-graph with one object. Thus we have shown
\begin{lemma}
Let $T$ be an accessible monad on $\ca GV$ over $\Set$ and $V$ be locally presentable. Suppose that for all $A$ and $B$ in $\ca G(V)^T$ one has comparisons \[ \kappa_{A,B}: A \BinFreeProd B \rightarrow A \times B \] in $\ca G(V)^T$. Then the unit $(T0,\mu_0)$ of the free product on $\ca G(V)^T$ is terminal.
\end{lemma}
\noindent and so we make
\begin{definition}
A monad $T$ on $\ca GV$ over $\Set$ is \emph{well-pointed} when $T0=1$.
\end{definition}
\begin{examples}\label{ex:well-pointed}
Let $E$ be a distributive multitensor on $V$ a category with coproducts and a terminal object. To say that the monad $\Gamma{E}$ is well-pointed is to say that the object $E_0$ of $V$ is terminal, by the explicit description of $\Gamma{E}$. Thus in particular if $T$ is a coproduct preserving monad on $V$, the multitensor $T^{\times}$ satisfies this property, and so $\Gamma{T^{\times}}$ is well-pointed. Thus the monads $\ca T_{{\leq}n}$ for strict $n$-categories for all $0{\leq}n{\leq}\infty$ are well-pointed.
\end{examples}
Whenever the unit of a symmetric monoidal category $(\ca W,\tensor)$ with finite products is terminal, one has natural comparisons
\[ \begin{array}{c} {\kappa_{(A_i)_i} : \Tensor\limits_i A_i \rightarrow \prod\limits_i A_i} \end{array} \]
defined by the condition that the composites $p_j\kappa_{(A_i)_i}$ are equal to the composites
\[ \xygraph{!{0;(2,0):} {\Tensor\limits_i A_i}="l" [r(1.2)] {\Tensor\limits_i 1|_jA_j}="m" [r] {A_j}="r" "l":@<1ex>"m"^-{\Tensor\limits_it_{A_i}|_j\id}:@{.>}@<.5ex>"r"} \]
for all $1{\leq}j{\leq}n$, where $p_j$ is the $j$-th product projection, the $t_{A_i}$ denote the unique maps into $1$ and dotted maps are the unit coherence isomorphisms in view of the fact that the unit of $\ca W$ is terminal. The compatibility between $\kappa$ and the monoidal structures on $\ca W$ involved is expressed by the following result, whose proof is a simple exercises in the definitions.
\begin{proposition}\label{prop:kappa-multitensor}
Let $(\ca W,\tensor)$ be a symmetric monoidal category with finite products whose unit is terminal. Then the components $\kappa_{(A_i)_i}$ form the coherences $(1_{\ca W},\kappa) : (\ca W,\prod) \rightarrow (\ca W,\tensor)$ of a symmetric lax monoidal functor.
\end{proposition}
\noindent In particular proposition(\ref{prop:kappa-multitensor}) implies that the $\kappa_{(A_i)_i}$ are the components of a morphism of multitensors $\Tensor \rightarrow \prod$, in fact this last statement \emph{is} the analogue of proposition(\ref{prop:kappa-multitensor}) in which the symmetries are disregarded. Instantiating to the case where $\tensor$ is the free product on $\ca G(V)^T$ for $T$ a good enough monad, one has in addition that the components $\kappa_{(X_i)_i}$ are identities on objects.
\begin{proposition}\label{prop:compare}
If $T$ is an accessible well-pointed monad on $\ca GV$ over $\Set$ and $V$ is locally presentable, then given $T$-algebras $((X_i,x_i))_i$, one may construct the comparison maps
\[ \begin{array}{c} {\kappa_{(X_i)_i} : \FreeProd\limits_i X_i \rightarrow \prod\limits_i X_i} \end{array} \]
of proposition(\ref{prop:kappa-multitensor}) in $\ca G(V)^T$ so that they are identities on objects.
\end{proposition}
\begin{proof}
Refer to the proof of theorem(\ref{thm:smcc-from-monad}). Note first that by the explicit construction of the free product on $\ca GV$, the $n$-ary tensor product preserves identity on objects maps. Thus in the coequaliser (\ref{eq:algtensor-def}) the morphism $\Tensor\limits_ix_i$ is the identity on objects. Since the coherences exhibiting $T$ as a symmetric monoidal monad live over $\Set$, the map $\mu{T(\phi)}$ in the coequaliser (\ref{eq:algtensor-def}) is also the identity on objects. Now from the transfinite construction of coequalisers in $\ca G(V)^T$ in terms of colimits in $\ca GV$ discussed in section(7) of \cite{EnHopII}, the underlying map in $\ca GV$ of the coequaliser $q_{x_i}$ is built by first taking a coequaliser down in $\ca GV$, then a series of successive pushouts and colimits of chains to construct a transfinite sequence down in $\ca GV$, and one proceeds until the length of the chain reaches the rank of $T$. But each stage of this process involves taking colimits of connected diagrams in $\ca GV$ involving only maps that are identities on objects, so at each stage one may take the colimit in $\ca GV$ to be the identity on objects. Thus in this way by a transfinite induction argument, one can indeed construct the coequaliser $q_{x_i}$ as being an identity on objects map. Thus the universal multimaps $\tilde{q}_{x_i}$ constructed in the proof of theorem(\ref{thm:smcc-from-monad}) are also identities on objects. For each $1{\leq}j{\leq}n$ the multimap corresponding to $p_j\kappa_{(X_i)_i}$, where $p_j$ is projection onto the $j$-th factor, has object map given by projection onto the $j$-th factor by the explicit description of $\kappa$. Since the explicitly constructed universal multimaps for $\ca F(\ca G(V)^T)$ have identity objects maps, it follows that the maps $p_j\kappa_{(X_i)_i}$ themselves have object maps given as projection onto the $j$-th factor, and so the result follows.
\end{proof}
%

\subsection{The pushout formula}\label{ssec:pushout-formula}
In this section we generalise the pushout formula for the funny tensor product, recalled in section(\ref{sec:funny}), to our setting. In the $\Cat$-case all the functors involved in the pushout formula are identities on objects, and so another place where one can locate this pushout is in a fibre of $\Cat{\rightarrow}\Set$. This point of view is the key for how to describe the general situation.

Suppose that $A$ is a distributive category. Suppose furthermore that $1 \in A$ is connected and regard $A$ has being over $\Set$ via the representable $A(1,-)$, which we recall has a left adjoint $(-) \cdot 1$ given by taking copowers with $1$. The connectedness of $1$ says exactly that $A(1,-)$ is coproduct preserving. In the discussion of section(\ref{sec:multicategories-of-graphs}) we denoted by $a_0$ the set $A(1,a)$, but we shall not do this here, preferring instead to use the notation $a_0$ for the object $A(1,a) \cdot 1$ of $A$. So for each element $x:1{\rightarrow}a$ one has a coproduct inclusion $c_x:1{\rightarrow}a_0$, and the component $\iota_a:a_0{\rightarrow}a$ of the counit of $(-){\cdot}1 \ladj A(1,-)$ is defined uniquely by $\iota_ac_x=x$ for all $x \in a$. The connectedness of $1$ ensures that $\iota_a$ is inverted by $A(1,-)$.

Given any sequence $(a_1,...,a_n)$ of objects of $A$ and $1{\leq}j{\leq}n$, one has by the distributivity of $A$ a canonical isomorphism
\begin{equation}\label{eq:dist-iso} \begin{array}{c} {\prod\limits_i a_{i,0}|_ja_j \iso \left( \prod\limits_{i{\neq}j} a_{i,0} \right) \cdot a_j} \end{array} \end{equation}
which by definition provides us with the following reformulation of the multimaps of $\ca F(A)$.
\begin{lemma}\label{lem:FA-multimap-reform}
Suppose that $A$ is a distributive category in which the terminal object $1$ is connected, and regard $A$ as over $\Set$ using the representable $A(1,-)$. Then to give a multimap
\[ f : (a_1,...,a_n) \rightarrow b \]
in $\ca F(A)$ is to give maps $f_j:\prod\limits_{i}a_{i,0}|_ja_j{\rightarrow}b$ for all $1{\leq}j{\leq}n$ such that
\[ \begin{array}{c} {f_j \left(\prod\limits_i \id|_j \iota_{a_j}\right) = f_k \left(\prod\limits_i \id|_k \iota_{a_k}\right)} \end{array} \]
for all $1{\leq}j,k{\leq}n$.
\end{lemma}
\begin{proof}
The object map of $f$ may be identified with the common composite map $f_j \left(\prod\limits_i \id|_j \iota_{a_j}\right)$ using the fact that $A(1,-)$ preserves coproducts, and using (\ref{eq:dist-iso}) the linear parts of $f$ may be identified with the $f_j$. Under these correspondences, the equation the $f_j$ must satisfy in the statement of this result amounts to the compatiblility between the underlying multifunction and the linear parts of $f$ required by the definition of $\ca F(A)$.
\end{proof}
\noindent Obtaining the pushout formula for the free product is simply a matter of applying this result to the case of $A = \ca G(V)^T$ for appropriate $T$, and using the fact that the free product of $T$-algebras is what classifies multimaps in $\ca F(\ca G(V)^T)$ by definition.
\begin{proposition}\label{prop:pushout-formula}
Let $T$ be a coproduct preserving, accessible and well-pointed monad on $\ca GV$ over $\Set$, and suppose that $V$ is locally presentable. Then for $T$-algebras $(X_1,...,X_n)$ their free product is the width-$n$ pushout of the diagram in $\ca G(V)^T$ whose maps are
\[ \begin{array}{c} {\prod\limits_i {\id}|_{j} i_{X_j} : \prod\limits_i X_{i,0} \rightarrow \prod\limits_i X_{i,0}|_j X_j} \end{array} \]
for $1{\leq}j{\leq}n$
\end{proposition}
\begin{proof}
Since $T$ is coproduct preserving and accessible, and $\ca GV$ is locally c-presentable by corollary(5.12) of \cite{EnHopII}, $\ca G(V)^T$ is locally c-presentable by theorem(5.6) of \cite{EnHopII}. Thus $\ca G(V)^T$ satsifies the hypotheses demanded of the category $A$ in lemma(\ref{lem:FA-multimap-reform}) enabling us to reformulate multimaps in $\ca F(\ca G(V)^T)$ as the appropriate cocones. But by theorem(\ref{thm:free-product-of-algebras}), $\ca F(\ca G(V)^T)$ is representable and the induced tensor product on $\ca G(V)^T$ is by definition the free product. The universal multimaps which define the free product of $T$-algebras correspond, via the reformulation of lemma(\ref{lem:FA-multimap-reform}), to the required width-$n$ pushout diagram.
\end{proof}
\begin{remark}
Using an argument similar to that given in proposition(\ref{prop:compare}) one can establish that the width-$n$ pushout diagram of proposition(\ref{prop:pushout-formula}) may be assumed to consist soley of identity on objects maps, and then one may recover proposition(\ref{prop:compare}) as a consequence of proposition(\ref{prop:pushout-formula}), though with the added hypotheses used by that proposition, by inducing the canonical comparison maps from the pushout in the obvious way.
\end{remark}
%

\section{Dropping multitensors}\label{sec:dropping}

\subsection{Overview}
Any multitensor $E$ on the category of algebras $V^T$ of some monad $(V,T)$ which is ``closed''
in the sense to be defined below, has been lifted, in the sense of \cite{EnHopII} section(7), from a multitensor down in $V$. This ``dropped''
multitensor $E$ is easy to describe in terms of the data given at the beginning. As a basic example one can start with the Gray tensor product of 2-categories, and then recapture the corresponding multitensor on $\ca G^2\Set$, and thus the monad on $\ca G^3\Set$ for Gray categories. The main result of this section, the multitensor dropping theorem, is proved by a similar argument to that used by Steve Lack in the proof of theorem 2 of \cite{LkMonFinMon}. It will then be used in section(\ref{sec:sesqui-algebras}) to explicitly describe the monads and operads whose algebras are categories enriched in the free products constructed in section(\ref{sec:free}).

\subsection{Closed multitensors}\label{ssec:closed-multitensors}
The multitensor dropping theorem applies to closed multitensors. We now discuss this notion.
\begin{definition}\label{def:closed-multitensor}
Let $V$ be a locally presentable category. A multitensor $(E,\iota,\sigma)$ is \emph{closed} when $E$ preserves colimits in each variable.
\end{definition}
\begin{example}
If $E$ is a genuine tensor product, then closedness in the sense of definition(\ref{def:closed-multitensor}) corresponds to closedness in the usual sense because for a locally presentable category $V$, an endofunctor $V{\rightarrow}V$ is cocontinuous iff it has a right adjoint.
\end{example}
\begin{lemma}\label{lem:cart-cocont}
Let $\phi:S{\rightarrow}T$ be a cartesian transformation between lra endofunctors on $V$ a locally cartesian closed category. If $T$ is cocontinuous then so is $S$.
\end{lemma}
\begin{proof}
Because of the cartesianness of $\phi$, $S$ is isomorphic to the composite
\[ \xygraph{{V}="p1" [r] {V/T1}="p2" [r(1.25)] {V/S1}="p3" [r] {V}="p4" "p1":"p2"^-{T_1}:"p3"^-{\phi_1^*}:"p4"^-{\textnormal{dom}}}. \]
$T_1$ is cocontinuous because $T$ is, $\phi_1^*$ is cocontinuous because $V$ is locally cartesian closed and $\textnormal{dom}$ is left adjoint to pulling back along the unique map $S1{\rightarrow}1$. 
\end{proof}
\begin{example}\label{ex:non-sig2}
A closed multitensor $E$ on $\Set$ is the same thing as a non-symmetric operad. For given such an $E$ and denoting by $E_n$ the set $\opE\limits_{1{\leq}i{\leq}n}1$, closedness gives
\[ \begin{array}{c} {\opE\limits_i X_i \iso \opE\limits_i \coprod\limits_{x \in X_i} 1 \iso E_n \times \prod\limits_i X_i.} \end{array} \]
Similarly the unit and substitution for $E$ determine the unit and substitution maps making the sequence $(E_n \, : \, n \in \N)$ of sets into an operad. Conversely given an operad $(E_n \, : \, n \in \N)$ of sets, the multitensor with object map
\[ \begin{array}{lcr} {(X_i)_i} & {\mapsto} & {E_n \times \prod\limits_i X_i} \end{array} \]
of example(2.6) of \cite{EnHopI} is clearly closed.
\end{example}
%

\subsection{Fibrewise Beck}\label{ssec:fibrewise-Beck}
We now present the key lemma that enables us to adapt Lack's proof \cite{LkMonFinMon} to our situation. Let
\begin{equation}\label{eq:adjunction} \xygraph{!{0;(2,0):} {V}="alg" [r] {W}="und"
"alg":@<-1ex>"und"_-{U}|{}="top":@<-1ex>"alg"_-{F}|{}="bot" "top":@{}"bot"|{\perp}} \end{equation}
be an adjunction, denote by $T$ the induced monad on $W$ and by $K:V{\rightarrow}W^T$ the functor induced by the universal property of $U^T:W^T{\rightarrow}W$. One form of the Beck theorem says: $U$ creates any reflexive coequalisers that it sends to absolute coequalisers iff $K$ is an isomorphism, and in this paper we call such right adjoints $U$ \emph{monadic}. Note however that elsewhere in the literature monadicity is often taken to mean that $K$ is an just equivalence of categories rather than an isomorphism.

Let us recall the relevant aspects of the proof of the Beck theorem (see \cite{Mac71} or any other textbook on category theory for a complete proof). First one analyses $U^T$ directly to see that it creates the appropriate coequalisers, and then conclude that when $K$ is an isomorphism $U$ obtains the desired creation property. For the converse, for a given $X \in V$ one notes that $U$ sends the diagram
\begin{equation}\label{eq:Beck-coeq} 
\xygraph{!{0;(2,0):} {FUFUX}="l" [r(1.5)] {FUX}="m" [r] {X}="r" "l":@<-2ex>"m"_-{FU\varepsilon_X}:"l"|-{F\eta_{UX}}:@<2ex>"m"^-{\varepsilon_{FUX}}:"r"^-{\varepsilon_X}} \end{equation}
to a split coequaliser, hence an absolute one, and then the creation property of $U$ can be used to deduce that $K$ is an isomorphism. If in a given situation the coequalisers of (\ref{eq:Beck-coeq}) satisfy some other useful and easily identifiable condition, then one immediately obtains a refinement of the Beck theorem in which one restricts attention to just those coequalisers.

The situation of interest to us is when the adjunction (\ref{eq:adjunction}) lives in $\CAT/\ca E$ for some fixed category $\ca E$, rather than just in $\CAT$. In this case for each $X \in \ca E$, the adjunction $F \ladj U$ restricts to an adjunction
\[ \xygraph{!{0;(2,0):} {V_X}="alg" [r] {W_X}="und" "alg":@<-1ex>"und"|{}="top"_-{U_X}:@<-1ex>"alg"|{}="bot"_-{F_X} "top":@{}"bot"|{\perp}} \]
where $V_X$ (resp, $W_X$) is the subcategory of $V$ (resp. $W$) consisting of the objects and arrows sent by the functor into $\ca E$ to $X$ and $1_X$. The induced monad $T_X$ on $W_X$ may be obtained by restricting $T$ in the same way. Note that any of the coequalisers (\ref{eq:Beck-coeq}) live in some $V_X$, because the components of $\eta$ all live in some $W_X$. Let us call any coequaliser in $V$ living in some $V_X$ a \emph{fibrewise coequaliser}. So any coequaliser of the form (\ref{eq:Beck-coeq}) is fibrewise. Thus one may restrict attention to such coequalisers in the proof of the Beck theorem, and so obtain the following fibrewise version of the Beck theorem.
\begin{lemma}\label{lem:fibrewise-Beck}
\emph{(Fibrewise Beck Theorem)}. Let $\ca E$ be a category. For a given adjunction $F \ladj U$ in $\CAT/\ca E$ TFSAE:
\begin{enumerate}
\item The functor $U$ is monadic in the usual sense.
\item For all $X \in \ca E$, $U_X$ is monadic.
\item $U$ creates any fibrewise reflexive coequaliser which it sends to an absolute coequaliser.
\end{enumerate}
\end{lemma}
%

\subsection{The multitensor dropping theorem}\label{ssec:dropping-theorem}
Recall that any multitensor $E$ on a category $V$ contains in particular a unary tensor product $E_1$. Rather than being trivial (and thus not mentioned) as in the case where $E$ is a genuine tensor product, $E_1$ is in general a monad on $V$. The central result of \cite{EnHopII}, theorem(7.3), explains how to ``lift'' a multitensor $E$ on a category $V$, to a multitensor $E'$ on the category $V^{E_1}$ whose unary part is trivial (though $E'$ could still be lax in general). Now we consider the reverse process -- given a monad playing the role of $E_1$, and a multitensor $E'$ on the category of algebras of $E_1$, we shall see how to recover $E$.

Suppose that $U:V{\rightarrow}W$ is monadic, accessible and coproduct preserving, with left adjoint denoted as $F$, $W$ is locally presentable, and $(E,u,\sigma)$ is a closed normal multitensor on $V$. Recall that by \cite{EnHopII} section(3.4) one has a multitensor $UEF$ on $W$ with object part $U\opE\limits_iFX_i$, and since $U$ preserves coproducts, $UEF$ is distributive. With $V$ the category of algebras of an accessible monad, it is locally presentable, and so by \cite{EnHopII} theorem(7.3) one obtains a multitensor $(UEF)'$ on $V \iso W^{(UEF)_1}$. Given that $E$ is normal it makes sense to ask whether $(UEF)' \iso E$. By the uniqueness part of \cite{EnHopII} theorem(7.3) this is the same as asking whether the composite functor
\[ \xygraph{{\Enrich E}="l" [r(1.25)] {\ca GV}="m" [r] {\ca GW}="r" "l":"m"^-{U^E}:"r"^-{\ca GU}} \]
is monadic. We now establish that this is indeed the case.
\begin{theorem}\label{thm:drop-mult}
Let $U:V{\rightarrow}W$ be monadic, accessible and coproduct preserving with left adjoint denoted as $F$, $W$ be locally presentable, and $(E,u,\sigma)$ be a closed normal multitensor on $V$. Then \[ (UEF)' \iso E \] as multitensors.
\end{theorem}
\begin{proof}
As argued above we must show that the composite $\ca G(U)U^E$ is monadic. Note that via the forgetful functors from $\ca GV$ and $\ca GW$ into $\Set$, all the monads and adjunctions involved in the present situation live in $\CAT/\Set$. Thus by lemma(\ref{lem:fibrewise-Beck}) it suffices to show that $\ca G(U)U^E$ creates any fibrewise reflexive coequalisers that it sends to absolute coequalisers. So we fix a set $X$ let
\[\xygraph{!{0;(1.5,0):} {A}="a" [r] {B}="b" "a":@<1.5ex>"b"^-{f}:"a"|{i}:@<-1.5ex>"b"_-{g}} \]
be a diagram of $E$-categories and $E$-functors in which the object maps of all the $E$-functors involved are $1_X$, and let
\[\xygraph{!{0;(2,0):} {\ca G(U)U^EA}="a" [r(1.5)] {\ca G(U)U^EB}="b" [r] {C}="c" "a":@<1ex>"b"^-{\ca G(U)U^Ef} "a":@<-1ex>"b"_-{\ca G(U)U^Eg}:"c"^-{h}} \]
be an absolute coequaliser of $W$-graphs. Since the object maps of $f$, $i$ and $g$ are identities, one may compute $C$ and $h$ as follows: take the object set of $C$ to be $X$ and the object map of $h$ to be the identity, and for all $a,b \in X$ take a coequaliser
\[\xygraph{!{0;(2,0):} {UA(a,b)}="a" [r] {UB(a,b)}="b" [r] {C(a,b)}="c" "a":@<1ex>"b"^-{Uf_{a,b}} "a":@<-1ex>"b"_-{Ug_{a,b}}:"c"^-{h_{a,b}}} \]
in $W$ to compute the hom $C(a,b)$ and the hom map $h_{a,b}$. We must exhibit a unique $h':B{\rightarrow}C'$ such that $\ca G(U)U^Eh'{=}h$ and
\[\xygraph{!{0;(2,0):} {A}="a" [r] {B}="b" [r] {C'}="c" "a":@<1ex>"b"^-{f} "a":@<-1ex>"b"_-{g}:"c"^-{h'}} \]
is a coequaliser of $E$-categories. The equation $\ca G(U)U^Eh'{=}h$ forces the object set of $C'$ to be $X$ and the object map of $h'$ to be the identity. By the monadicity of $U$, one induces $h'_{a,b}:B(a,b){\rightarrow}C'(a,b)$ as the unique map in $V$ coequalising $f_{a,b}$ and $g_{a,b}$. So far we have constructed the underlying $V$-graph of $C'$, which we shall also denote as $C'$, and the underlying $V$-graph morphism of $h'$ which we shall also denote by $h'$. In fact by the uniqueness part of the monadicity of $U$, the $V$-graph morphism $h'$ is forced to be as we have just constructed it. Thus to finish the proof it suffices to do two things: (1) show that there is a unique $E$-category structure on $C'$ making $h'$ an $E$-functor, and (2) show that $h'$ is indeed the coequaliser $f$ and $g$ in $\Enrich E$.

Let us now witness the unique $E$-category structure. Let $n \in \N$ and $x_0,...,x_n$ be elements of $X$. Then one induces the corresponding composition map for $C'$ from those of $A$ and $B$ as shown
\begin{equation}\label{eq:def-Cpr} \xygraph{!{0;(3,0):(0,.5)::} {\opE\limits_i A(x_{i{-}1},x_i)}="a1" [r] {\opE\limits_i B(x_{i{-}1},x_i)}="b1" [r] {\opE\limits_i C'(x_{i{-}1},x_i)}="c1" [d] {C'(x_0,x_n)}="c2" [l] {B(x_0,x_n)}="b2" [l] {A(x_0,x_n)}="a2"
"a1":@<1ex>"b1"^-{\opE\limits_i f_{x_{i{-}1},x_i}} "a1":@<-1ex>"b1"_-{\opE\limits_i g_{x_{i{-}1},x_i}}:"c1"^-{\opE\limits_i h'_{x_{i{-}1},x_i}} "a2":@<1ex>"b2"^-{f_{x_0,x_n}} "a2":@<-1ex>"b2"_-{g_{x_0,x_n}}:"c2"^-{h'_{x_0,x_n}} "a1":"a2" "b1":"b2" "c1":@{.>}"c2"} \end{equation}
because the top row is a coequaliser by the closedness of $E$ and the $3 \times 3$-lemma (see corollary(7.11) of \cite{EnHopII}). Note that the putative $E$-category structure is by definition uniquely determined by the condition that $h'$ becomes an $E$-functor. That is given sequences $(x_{ij} \,\, : \, \, 1{\leq}j{\leq}n_i)$ for each $1{\leq}i{\leq}n$ such that $x_{i0}=x_{i-1}$ and $x_{in_i}=x_i$ from $X$, we must verify the commutativity of the corresponding
\begin{equation}\label{eq:Cpr-ax} \xygraph{!{0;(3,0):(0,.333)::} {\opE\limits_i\opE\limits_j C'(x_{(ij){-}1},x_{ij})}="tl" [r] {\opE\limits_{ij} C'(x_{(ij){-}1},x_{ij})}="tr" [d] {C'(x_0,x_n)}="br" [l] {\opE\limits_{i}C'(x_{i{-}1},x_i)}="bl" "tl"(:@<1ex>"tr"^-{\sigma}:"br",:"bl":"br")} \end{equation}
Here is a thumbnail sketch of the diagram which enables one to witness this:
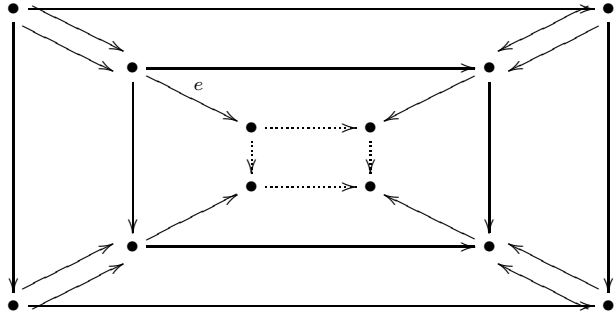
\[ \xygraph{!{0;(1.25,0):(0,.5)::} {\bullet}="otl" [dr] {\bullet}="mtl" [dr] {\bullet}="itl" [r] {\bullet}="itr" ([ur] {\bullet}="mtr" [ur] {\bullet}="otr",[d] {\bullet}="ibr" ([dr] {\bullet}="mbr" [dr] {\bullet}="obr",[l] {\bullet}="ibl" [dl] {\bullet}="mbl" [dl] {\bullet}="obl"))
"otl":@<1ex>"mtl" "otl":@<-1ex>"mtl":"itl"^-{e} "otr":@<1ex>"mtr" "otr":@<-1ex>"mtr":"itr" "obl":@<1ex>"mbl" "obl":@<-1ex>"mbl":"ibl" "obr":@<1ex>"mbr" "obr":@<-1ex>"mbr":"ibr" "otl"(:"otr":"obr",:"obl":"obr") "mtl"(:"mtr":"mbr",:"mbl":"mbr") "itl"(:@{.>}"itr":@{.>}"ibr",:@{.>}"ibl":@{.>}"ibr")} \]
In this diagram the inner rectangle is (\ref{eq:Cpr-ax}). The outer and middle rectangles are the corresponding axioms for $A$ and $B$ respectively and so commute by definition. The diagonal parts of the diagram are coequalisers by the closedness of $E$ and the $3 \times 3$-lemma. The bottom left and bottom right coequalisers in this diagram coincide with the top and bottom row of (\ref{eq:def-Cpr}) respectively. The vertical and horizontal maps provide natural transformations of coequaliser diagrams by definition, modulo the commutativity of the inner rectangle to be established. But by the commutativities just witnessed, this rectangle does commute after precomposition with $e$, but as a coequalising map $e$ is an epimorphism, and so it does indeed commute. This concludes the verification of (1).

As for (2) suppose that one has an $E$-functor $k:B{\rightarrow}D$ such that $kf=kg$. We must exhibit a unique $E$-functor $k':C'{\rightarrow}D$ such that $k'h'=k$. The underlying $V$-graph map is determined uniquely with this condition because $U^Eh'$ is the coequaliser of $U^Ef$ and $U^Eg$ by construction, so to finish the proof it suffices to verify that this $V$-graph map $k'$ is compatible with the $E$-category structures. But for a given $(x_0,...x_n)$ as above, by definition the corresponding $E$-functoriality rectangle does commute after precomposition with the map $\opE\limits_ih'_{x_{i{-}1},x_i}$, and so the result follows since this map is an epimorphism.
\end{proof}
%

\section{Monads and operads for sesqui-algebras}\label{sec:sesqui-algebras}

\subsection{Sesqui-algebras}\label{ssec:sesqui-algebras} Given a monad $T$ on $\ca GV$ over $\Set$ satisfying the hypotheses of theorem(\ref{thm:free-product-of-algebras}), as we have seen one may speak about the free product of $T$-algebras. One may then consider categories enriched in $\ca G(V)^T$ via the free product. We shall call such structures, that is categories enriched in $\ca G(V)^T$ for the free product, \emph{sesqui-$T$-algebras}. When $T$ is the category monad on $\ca G\Set$, sesqui-$T$-algebras are just sesqui-categories in the usual sense.

\subsection{Monads for sesqui-algebras}\label{ssec:monads-sesqui-algebras} 
Since $T$ is a monoidal monad with respect to the free product on $\ca GV$ one has a multitensor on $\ca GV$ with object maps
\[ \begin{array}{c} {(X_i)_i \mapsto T\FreeProd\limits_iX_i} \end{array} \]
as discussed in \cite{EnHopII} section(6.6). However by remark(\ref{rem:tensor-of-frees}) $F^T$ is strict monoidal with respect to the free products on $\ca GV$ and $\ca G(V)^T$ and so one has an equality
\[ \begin{array}{c} {T\FreeProd = U^T\FreeProd{F^T}} \end{array} \]
of multitensors. From the right hand side of this last equation we see that we are in the situation of the multitensor dropping theorem of section(\ref{ssec:dropping-theorem}) if our monad $T$ is good enough.
\begin{theorem}\label{thm:monad-for-sesqui-algebras}
Suppose that $T$ is an accessible and coproduct preserving monad on $\ca GV$ over $\Set$ and let $V$ be locally presentable.
\begin{enumerate}
\item  The free product of $T$-algebras may be recovered as the lifted multitensor $(T\FreeProd)'$.\label{free-as-lift}
\item  The monad on $\ca G^2V$ whose algebras are sesqui-$T$-algebras is given explicitly as $\Gamma(T\FreeProd)$.\label{sesqui-monad}
\end{enumerate}
\end{theorem}
\begin{proof}
(\ref{free-as-lift}) is immediate from theorem(\ref{thm:free-product-of-algebras}) and theorem(\ref{thm:drop-mult}), and so (\ref{sesqui-monad}) follows by \cite{EnHopII} corollary(4.9).
\end{proof}
\noindent So we have obtained in (\ref{sesqui-monad}) an explicit combinatorial description of the monad on $\ca G^2V$ whose algebras are sesqui-$T$-algebras.

\subsection{Operads for sesqui-algebras}\label{ssec:operads-sesqui-algebras}
Applying theorem(\ref{thm:monad-for-sesqui-algebras}) to the monad $T=\ca T_{{\leq}1}$ for categories, one obtains an explicit description of the monad on 2-globular sets whose algebras are sesqui-categories. In this case we know a bit more, namely that this monad is part of a 2-operad, in fact it is one of the basic examples of such -- see \cite{Bat98}.

More generally given an $n$-operad $A{\rightarrow}\ca T_{{\leq}n}$, it is intuitively obvious that sesqui-$A$-categories are describable by an $(n{+}1)$-operad, because such a structure on an $(n{+}1)$-globular set amounts to a category structure on the underlying graph and $A$-algebra structures on the homs, with no compatibility between them. In other words the structure of a sesqui-$A$-algebra may be described by reinterpretting the compositions and axioms for $A$-algebras one dimension higher, giving a category structure in dimensions $0$ and $1$ and imposing no further axioms, and so the data and axioms are inherently of the type describable by an $(n{+}1)$-operad.

However the explanation just given is not really a rigorous proof that sesqui-$A$-algebras are $(n{+}1)$-operadic. To give one it is necessary to provide a cartesian monad morphism as on the left, or equivalently a cartesian multitensor map as on the right in
\[ \begin{array}{lccr} {\Gamma(A\FreeProd) \rightarrow \ca T_{{\leq}{n{+}1}}} &&& {A\FreeProd \rightarrow \ca T^{\times}_{{\leq}n}.} \end{array} \]
In this section we shall give the general construction, the setting being a good enough monad $T$ on $\ca GV$ over $\Set$ replacing $\ca T_{{\leq}n}$. As we shall see, one of the fundamental properties required of $T$ is that it be well-pointed.

Let $V$ be locally presentable and $T$ be a monad on $\ca GV$ over $\Set$ which is accessible and well-pointed. As justified by remark(\ref{rem:tensor-of-frees}) we regard $F^T$ as a strict monoidal functor $(\ca GV,\FreeProd){\rightarrow}(\ca G(V)^T,\FreeProd)$. So for any sequence of objects $(X_i)_i$ of $\ca GV$, we have
\[ T\FreeProd\limits_iX_i = U^T\FreeProd\limits_i F^TX_i \]
as we saw in section(\ref{ssec:sesqui-algebras}). From section(\ref{ssec:free-vs-cartesian}) we have the comparison
\[ \begin{array}{c} {\kappa_{(F^TX_i)_i} : \FreeProd\limits_i F^TX_i \rightarrow \prod\limits_i F^TX_i} \end{array} \]
in $\ca G(V)^T$. We define $\kappa_{(X_i)_i} = U^T\kappa_{(F^TX_i)_i}$, and since $U^T$ preserves products, we regard it as a map
\[ \begin{array}{c} {\kappa_{(X_i)_i} : T\FreeProd\limits_i X_i \rightarrow \prod\limits_i TX_i} \end{array} \]
in $\ca GV$. By definition $\kappa$ is natural in the $X_i$. The key lemma of this section is
\begin{lemma}\label{lem:key}
Let $V$ be locally presentable and extensive, and $T$ be a monad on $\ca GV$ over $\Set$ which is accessible, well-pointed, lra, distributive and path-like. Then the maps $\kappa_{(X_i)_i}$ are cartesian natural in the $X_i$.
\end{lemma}
\noindent and we defer the proof of this result until after we have discussed its consequences. First note that lemma(\ref{lem:key}) exhibits $T\FreeProd$ as a $T$-multitensor by proposition(\ref{prop:compare}). Most importantly from lemma(\ref{lem:key}) the central result of this section follows immediately.
\begin{theorem}\label{thm:sesqui-operadic}
Let $V$ be locally presentable and extensive, and $T$ be a monad on $\ca GV$ over $\Set$ which is accessible, well-pointed, lra, distributive and path-like. Suppose that $\psi:A{\rightarrow}T$ is a $T$-operad. Then
\[ \xygraph{!{0;(2.5,0):} {A\FreeProd\limits_iX_i}="l" [r] {T\FreeProd\limits_iX_i}="m" [r] {\prod\limits_iTX_i}="r" "l":@<1ex>"m"^-{\psi_{\FreeProd\limits_iX_i}}:@<1ex>"r"^-{\kappa_{(X_i)_i}}} \]
are the components of a $T$-multitensor.
\end{theorem}
\begin{proof}
By proposition(\ref{prop:kappa-multitensor}) and since $\psi$ is a morphism of monads, the given composite maps form a morphism of multitensors. Its cartesianness follows since $\psi$ is cartesian by definition and $\kappa_{(X_i)_i}$ is cartesian natural in the $X_i$ by lemma(\ref{lem:key}).
\end{proof}
\noindent Applying this result in the cases $T=\ca T_{{\leq}n}$ we see that for an $n$-operad $A$, sesqui-$A$-algebras are indeed describable by an $(n{+}1)$-operad. The rest of this section is devoted to proving lemma(\ref{lem:key}).

Let us denote by $\tilde{\kappa}_{(X_i)_i}$ the composite
\[ \xygraph{!{0;(2.5,0):} {\FreeProd\limits_iX_i}="l" [r] {T\FreeProd\limits_iX_i}="m" [r] {\prod\limits_iTX_i}="r" "l":@<1ex>"m"^-{\eta_{\FreeProd\limits_iX_i}}:@<1ex>"r"^-{\kappa_{(X_i)_i}}} \]
which by definition is also natural in the $X_i$. It turns out that for lemma(\ref{lem:key}) it suffices to consider the cartesian naturality of $\tilde{\kappa}$.
\begin{lemma}\label{lem:tilde-reduction}
Let $V$ and $T$ satisfy the hypotheses of lemma(\ref{lem:key}). Then $\kappa$ is cartesian-natural in the $X_i$ iff $\tilde{\kappa}$ is.
\end{lemma}
\begin{proof}
If $\kappa$ is cartesian natural then so is $\tilde{\kappa}$ since $\eta$ is cartesian. For the converse note that by definition the $\kappa_{(X_i)_i}$ are $T$-algebra maps, and thus one may recover $\kappa_{(X_i)_i}$ as the composite
\[ \xygraph{!{0;(2,0):} {T\FreeProd\limits_i X_i}="p1" [r] {T\prod\limits_iTX_i}="p2" [r] {\prod\limits_iT^2X_i}="p3" [r] {\prod\limits_iTX_i}="p4" "p1":@<1ex>"p2"^-{T\tilde{\kappa}}:@<1ex>"p3"^-{k_{TX_i}}:@<1ex>"p4"^-{\prod\limits_i \mu_{X_i}}} \]
where $k_{TX_i}$ is the product obstruction map for $T$. This is so since the composite of the last two arrows in this string is the $T$-algebra structure of $\prod\limits_iTX_i$. By \cite{Fam2fun} lemma(2.15) $k_{TX_i}$ is cartesian natural in the $X_i$ since $T$ is lra. Since $\mu$ is cartesian and a product of pullback squares is a pullback square, each of the maps in the above composite is cartesian natural in the $X_i$, and so the result follows.
\end{proof}
The components $\tilde{\kappa}_{(X_i)_i}$ are identities on objects maps, and so it suffices by \cite{EnHopII} lemma(5.17) to show that in any naturality square, all the induced commutative squares on the homs are pullbacks in $V$.

Recall that the functor $(-)_0:\ca GV{\rightarrow}\Set$ which sends a $V$-graph to its set of objects has a representing object $0$, and so we have an adjunction
\[ \xygraph{!{0;(2,0):} {\ca GV}="l" [r] {\Set}="r" "l":@<-1ex>"r"_-{(-)_0}|-{}="b":@<-1ex>"l"_-{(-){\cdot}0}|-{}="t" "t":@{}"b"|{\perp}} \]
and we regard the $V$-graph $Z \cdot 0$ as having object set $Z$ and all homs equal to $\emptyset$ the initial object of $V$. For convenience we write the functor $(-) \cdot 0$ as though it were an inclusion. The counit of this adjunction has components we will denote as $i_X:X_0 \rightarrow X$, and these are the identity on objects, and the hom maps are determined uniquely. Obviously $i_X$ is cartesian natural in $X$: to see this it suffices by \cite{EnHopII} lemma(5.17) to look at all the induced squares between homs of a given naturality square since $i_X$ is the identity on objects, and these are all of the form
\[ \xygraph{{\emptyset}="tl" [r] {A}="tr" [d] {B}="br" [l] {\emptyset}="bl" "tl"(:"bl":"br",:"tr":"br")} \]
which since $\emptyset$ is a strict initial object by the extensivity of $V$, are automatically pullbacks.

The extensivity of $V$ enables a more efficient description of the free product on $\ca GV$. Recall that the general definition of the homs of $\FreeProd\limits_iX_i$ required a case split -- see the formula(\ref{GV-tensor}) of section(\ref{ssec:fp-GV}). Let us consider the functor
\[ d : \Set \rightarrow \ca GV \]
which sends a set $Z$ to the $V$-graph with object set $Z$ and homs defined by $dZ(a,b)=\emptyset$ if $a \neq b$ and $dZ(a,b) = 1$ if $a = b$. Clearly $d$ may be identified with taking copowers with $1$, that is $dZ \iso Z \cdot 1$. For convenience we shall also denote by $d$ the endofunctor $d(-)_0$ of $\ca GV$, in other words for $X \in \ca GV$ we write $dX$ for $d(X_0)$. Since the cartesian product of $V$ is distributive, the formula(\ref{GV-tensor}) of section(\ref{ssec:fp-GV}) may be re-expressed as in
\begin{lemma}
Let $V$ have products and coproducts, with coproducts distributing over finite products. Then the homs of the free product on $\ca GV$ may be re-expressed as
\begin{equation}\label{GV-tensor2} \begin{array}{c} {\left(\FreeProd\limits_{1{\leq}i{\leq}n} X_i\right)(a,b) = \coprod\limits_{1{\leq}j{\leq}n}\prod\limits_{1{\leq}i{\leq}n} dX_i(a_i,b_i)|_{j}X_{j}(a_j,b_j)} \end{array} \end{equation}
\end{lemma}
\begin{proof}
If $a=b$ then the homs $dX_i(a_i,b_i)$ are all $1$ making the right hand side the coproduct of the $X_j(a_j,b_j)$. If all but one, say the $j$-th, coordinates of $a$ and $b$ coincide, then the only non-$\emptyset$ summand on the right hand side of (\ref{GV-tensor2}) is $X_j(a_j,b_j)$. Otherwise all summands are $\emptyset$. Thus (\ref{GV-tensor2}) coincides with formula(\ref{GV-tensor}) of section(\ref{ssec:fp-GV}) under the given hypotheses.
\end{proof}
Assuming that $T$ is well-pointed and preserves coproducts we have
\[ dX \iso X_0 \cdot 1 \iso X_0 \cdot T0 \iso T(X_0 \cdot 0) \]
and so $dX$ is in fact a free $T$-algebra. For convenience we write $dX = TX_0$. Note that we have maps
\[ Ti_X : dX \rightarrow TX \]
which given the extensivity of $V$ and lra'ness of $T$, are cartesian-natural in $X$. Thus we can define maps
\[ \begin{array}{c} {\overline{\kappa}_{(X_i)_i} : \FreeProd\limits_i X_i \rightarrow \prod\limits_i TX_i} \end{array} \]
as follows. They are identities on objects. To define the hom maps it suffices in view of (\ref{GV-tensor2}) to define the maps
\[ \begin{array}{c} {\left(\overline{\kappa}_{(X_i)_i}\right)_{a,b} c_j : \prod\limits_i dX_i(a_i,b_i)|_{j}X_{j}(a_j,b_j) \rightarrow \prod\limits_i TX_i(a_i,b_i)} \end{array} \]
where $c_j$ denotes the $j$-th coproduct inclusion, and we define
\[ \begin{array}{c} {\left(\overline{\kappa}_{(X_i)_i}\right)_{a,b} c_j = \prod\limits_i (Ti_X)_{a_i,b_i}|_j (\eta_{X_j})_{a_j,b_j}} \end{array} \]
and by this definition $\overline{\kappa}$ is clearly natural in its arguments. In fact
\begin{lemma}\label{lem:kappa-bar}
Let $V$ be complete and extensive, and suppose that the monad $T$ on $\ca GV$ over $\Set$ is well-pointed, cartesian and preserves coproducts. Then $\overline{\kappa}_{(X_i)_i}$ is cartesian natural in the $X_i$.
\end{lemma}
\begin{proof}
By extensivity it suffices to show that the $\left(\overline{\kappa}_{(X_i)_i}\right)_{a,b} c_j$ are cartesian natural, and this follows since $i_X$ and $\eta_X$ are cartesian, and products of pullback squares are pullbacks.
\end{proof}
To conclude the proof of lemma(\ref{lem:key}) it remains to show that $\overline{\kappa}=\tilde{\kappa}$. The problem with doing this at the generality of lemma(\ref{lem:kappa-bar}) is that $\overline{\kappa}$ and $\tilde{\kappa}$ implicitly involve the two different ways of describing multimaps of $V$-graphs -- the general one from $\ca F(\ca GV)$ involving object maps and linear parts, versus the description of example(\ref{ex:FGV}) involving object maps and individual hom maps. Without the more explicit description of $T$ that becomes available since it is distributive and path-like, it doesn't seem possible to identify $\overline{\kappa}$ and $\tilde{\kappa}$.

So we now use that $T$ is distributive and path-like so that there exists a distributive multitensor $E$ on $V$ and $\Gamma{E} \iso T$ by proposition(4.11). This means that for all $X \in \ca GV$ and $a,b \in X$ we have maps
\[ c_{(x_i)_i} : \opE\limits X(x_{i{-}1},x_i) \rightarrow TX(a,b) \]
for all sequences $(x_i)_i$ starting at $a$ and finishing at $b$, and taken together these form a coproduct cocone. Thus given a multimap
\[ f : (X_1,...,X_n) \rightarrow Y \]
of $V$-graphs, the multimap $Tf : (TX_i)_i{\rightarrow}TY$ may be described explicitly as having object map that of $f$, and hom map corresponding to $x \in(X_i)_{i{\neq}i^*}$ and $a,b \in X_{i^*}$, defined by the commutativity of
\[ \xygraph{!{0;(0,1):(0,4)::} {\opE\limits_jX_{i^*}(z_{j{-}1},z_j)}="tl" [r] {TX_{i^*}(a,b)}="tr" [d] {TY(f(x|_{i^*}a),f(x|_{i^*}b))}="br" [l] {\opE\limits_jY(f(x|_{i^*}z_{j{-}1}),f(x|_{i^*}z_j))}="bl" "tl"(:"tr"^-{c_{(z_j)_j}}:"br"^-{(Tf)_{x,a,b}},:@<1ex>"bl"_-{\opE\limits_jf_{x,z_{j{-}1},z_j}}:"br"_-{c_{(f(x|_{i^*}z_j))_j}})} \]
for all sequences $(z_j)_j$ in $X_{i^*}$ from $a$ to $b$. So using this extra explicit information we finish the proof of lemma(\ref{lem:key}) in
\begin{lemma}\label{lem:last}
Let $V$ and $T$ satisfy the hypotheses of lemma(\ref{lem:key}). Then $\tilde{\kappa}=\overline{\kappa}$.
\end{lemma}
\begin{proof}
Given $V$-graphs $(X_1,...,X_n)$ and $1{\leq}j{\leq}n$ we must show $p_j\tilde{\kappa}_{(X_i)_i}=p_j\overline{\kappa}_{(X_i)_i}$, where $p_j$ is the $j$-th projection of the product. We will show that the corresponding multimaps $(X_i)_i{\rightarrow}TX_j$ in $\ca F(\ca GV)$ are the same. To this end consider
\[ \xygraph{!{0;(2.5,0):(0,.5)::} {(X_i)_i}="tl" [r] {(TX_i)_i}="tm" [r] {(T0|_jTX_j)_i}="tr" [d] {TX_j}="br" [l] {T\FreeProd\limits_iX_i}="bm" [l] {\FreeProd\limits_iX_i}="bl"
"tl"(:"tm"^-{(\eta_{X_i})_i}(:"tr"^-{(t_{TX_i}|_j\id)_i}:"br"^{T(\alpha_{(0|_{j}X_j)_i})},:"bm"^{T\alpha_{(X_i)_i}}),:"bl"_{\alpha_{(X_i)_i}}:@<1ex>"bm"_-{\eta}:@<1ex>"br"_-{p_j\kappa_{(X_i)_i}})} \]
in which the linear maps are the unique maps in view of $T0 \iso 1$ by well-pointedness, and the $\alpha$ maps are the universal maps described just before proposition(\ref{prop:FGV-rep}). Note that the right hand square may be regarded as living in $\ca F(\ca G(V)^T)$ and it commutes by the definition of $\kappa$. The left hand square commutes by the multinaturality of $\eta$. By the explicit description of the hom maps of $T(\alpha_{(0|_{j}X_j)_i})$ in terms of the corresponding multitensor $E$ (see the discussion just before the statement of this lemma), the top composite of the diagram
\[ T(\alpha_{(0|_{j}X_j)_i})(t_{TX_i}\eta_{X_i}|_j\eta_{X_j})_i \]
is the multimap corresponding to $p_j\overline{\kappa}_{(X_i)_i}$, and so the outside of the diagram witnesses $p_j\tilde{\kappa}_{(X_i)_i}=p_j\overline{\kappa}_{(X_i)_i}$ at the level of multimaps.
\end{proof}

\section{Acknowledgements}\label{sec:Acknowledgements}
The author would like to acknowledge Michael Batanin, Clemens Berger, John Bourke, Albert Burroni, Denis-Charles Cisinski, Richard Garner, Tom Hirschowitz, Paul-Andr\'{e} Melli\`{e}s, Fran\c{c}ois M\'{e}tayer and Jacques Penon for interesting discussions on the substance of this paper. He would also like to express his gratitude to the laboratory PPS (Preuves Programmes Syst\`{e}mes) in Paris and the Max Planck Institute in Bonn for the excellent working conditions he enjoyed during this project.

\end{document}